\documentclass[a4paper,twocolumn]{article}
\pdfoutput=1
\usepackage{graphicx,amsmath,bm,mathtools,float,titlesec,authblk}
\usepackage[a4paper, margin=1.5cm]{geometry}
\usepackage{lmodern}
\usepackage[T1]{fontenc}

\titleformat*{\section}{\Large\bfseries}
\titleformat*{\subsection}{\large\bfseries}

\title{\textbf{Tractable nonlinear memory functions as a tool to capture and explain dynamical behaviours}}
\author[1,2]{Edgar Herrera-Delgado}
\author[1]{James Briscoe*}
\author[2,3]{Peter Sollich*}
\affil[1]{The Francis Crick Institute, 1 Midland Rd, London NW1 1AT, UK}
\affil[2]{Department of Mathematics, King's College London, Strand, London WC2R 2LS, UK}
\affil[3]{Institut f\"ur Theoretische Physik, University of G\"ottingen, Friedrich-Hund-Platz 1, 37077 G\"ottingen, Germany}
\date{}

\begin{document}

\newcommand{\xs}{\bm{x}^{\textrm{s}}}
\newcommand{\hxs}{\hat{\bm{x}}^{\textrm{s}}}
\newcommand{\xb}[1]{\bm{x}^{\textrm{b}#1}}
\newcommand{\fbs}[2]{f_{#1 s}(\xs, #2)}
\newcommand{\lbb}{l_{bb'}}
\newcommand{\expl}[2]{e^{\int_{#1}^{#2}d\tau\,\bm{l}(\phio)}}
\newcommand{\exps}[3]{E_{#1}({#2},{#3})}
\newcommand{\fo}[1]{f^0_{#1}(\bm{\phi}(\xs,\tau))}
\newcommand{\phio}{\bm{\phi}(\xs,\tau)}
\newcommand{\rc}[2]{\bm{R}^{\textrm{#1}}#2}
\newcommand{\F}{\bm{F}}
\newcommand{\A}{\bm{A}}
\newcommand{\zmv}{ZMn }
\newcommand{\zms}{ZMs }

\twocolumn[
  \begin{@twocolumnfalse}
    \maketitle
    \vspace{-0.8cm}
    james.briscoe@crick.ac.uk\\
	peter.sollich@uni-goettingen.de
    \begin{abstract}
	Mathematical approaches from dynamical systems theory are used in a range of fields. This includes biology where they are used to describe processes such as protein-protein interaction and gene regulatory networks. As such networks increase in size and complexity, detailed dynamical models become cumbersome, making them difficult to explore and decipher. This necessitates the application of simplifying and coarse graining techniques in order to derive explanatory insight. Here we demonstrate that Zwanzig-Mori projection methods can be used to arbitrarily reduce the dimensionality of dynamical networks while retaining their dynamical properties. We show that a systematic expansion around the quasi-steady state approximation allows an explicit solution for memory functions without prior knowledge of the dynamics. The approach not only preserves the same steady states but also replicates the transients of the original system. The method also correctly predicts the dynamics of multistable systems as well as networks producing sustained and damped oscillations. Applying the approach to a gene regulatory network from the vertebrate neural tube, a well characterised developmental transcriptional network, identifies features of the regulatory network responsible for its characteristic transient behaviour. Taken together, our analysis shows that this method is broadly applicable to multistable dynamical systems and offers a powerful and efficient approach for understanding their behaviour.
	\end{abstract}
	\vspace{0.4cm}
  \end{@twocolumnfalse}
  ]

\section*{Introduction}
In complex dynamical systems, comprising multiple interacting components, it can be difficult to identify causal mechanisms and to dissect the function of parts of a system. Nonlinearities and feedback complicate intuitive understanding and these difficulties increase with the size and complexity of a system. Examples include biological processes such as protein-protein interaction and gene regulatory networks
%are examples of such systems
\cite{Davidson2010,Snider2015}.
Mathematical models of these systems allow exploratory analysis and can provide insight but become less practical as system size grows. More importantly, the complexity can obscure the explanation for unexpected or emergent behaviours that originate in the dynamics of a system. For these reasons, a variety of approaches have been developed to reduce the complexity of models while preserving desired features of their behaviour. An important class of tools are dimensionality reduction techniques that coarse grain parts of a system \cite{Rega2005,Schnoerr2017,Bronstein2018}.

The Zwanzig-Mori formalism provides an exact dimensionality reduction of a dynamical system based on a separation %system
into an arbitrary ``subnetwork'', the components of which are tracked explicitly, and a ``bulk'' containing the components that are replaced with ``memory functions'' \cite{Nakajima1958,Zwanzig1961,Mori1965,Kawasaki1973}. These functions describe how the current subnetwork state feeds back, through the activity of molecular species in the bulk, to affect the subnetwork at a later time. This approach, specifically its nonlinear version, was originally developed for the dynamics of physical systems \cite{Zwanzig2001}, but later generalized by Chorin and coworkers \cite{Chorin2000,Chorin2002}, with related uses also in optimal coarse graining \cite{Weinan2008}. A limitation, however, is that the memory functions are generally impossible to calculate in closed form \cite{Chorin2000,Chorin2002}. Although approximate expressions can be derived in special cases \cite{Chorin2006a,Stinis2006,Beck2009,Thomas2012,Gouasmi2017}, this restricts the applicability of the formalism. An alternative is to map the nonlinear system to a physical system consistent with the original; while this can be effective it does not necessarily simplify the problem \cite{Xing2011,Mukhopadhyay2018}. Another option is to expand the dynamical equations around a fixed point and derive memory functions from this approximation \cite{Rubin2014,Bravi2017a}, but for multistable or oscillatory systems the memory functions obtained in this way do not capture all the behaviours of a system.

To address this limitation, we develop a method, based on the formalism of \cite{Chorin2000}, that allows the calculation of memory functions for generic dynamical systems without prior knowledge of the dynamics. We make one assumption: the bulk must not generate fixed points beyond those of the subnetwork, more specifically it must have a unique steady state for any subnetwork state, such as in \cite{Gouasmi2017}. This is a natural condition: the subnetwork must be able to produce all fixed points itself, otherwise coarse graining cannot succeed. As the starting point for the dimensionality reduced dynamics we use the Quasi-steady state (QSS) approximation, where the bulk is always in steady state with the current subnetwork state as has been used in other contexts \cite{Kang2019}. Memory functions are then constructed to correct the projected subnetwork state, by accounting for departures of the bulk from its steady state. Our main technical result is an explicit solution for the functions capturing these memory effects, derived in a systematic expansion around the QSS approximation. We demonstrate that the approach accurately predicts the dynamics of systems that produce multiple steady states and even sustained or damped oscillations. We also illustrate its use by applying it to a gene regulatory network from the embryonic vertebrate neural tube \cite{Cohen2014}. This is a transcriptional network of four interacting transcription factors with well described transient dynamics. We show how the memory functions generated by this approach provide insight into the features of the regulatory network that produce this transient behaviour. Taken together, the analysis introduces a broadly applicable method for the investigation and analysis of complex dynamical systems.

\section*{Mathematical derivation \label{sec:mathCh}}
\subsection*{Initial definitions}
Following \cite{Chorin2000}, we start from a system with degrees of freedom $\bm{x}$ evolving deterministically in time according to some nonlinear functions $\bm{R}$:
\begin{equation}
\frac{d\bm{x}}{dt}=\rc{}{(\bm{x})}\label{reaqsEx}
\end{equation}
We define the ``flow'' $\bm{\phi}(\bm{x},t)$ as the state the system reaches at time $t$ if it starts in some initial state $\bm{x}$; this function thus obeys $\bm{\phi}(\bm{x},0)=\bm{x}$ and $\frac{\partial}{\partial t}\bm{\phi}(\bm{x},t)=\rc{}(\bm{\phi}(\bm{x},t))$.  We want to understand the dynamics of some chosen set of observables that we denote by the vector $\bm{A}$. Such observables are functions of the state of the system, which we write as $\A(\bm{x})$. By analogy with the definition of $\bm{\phi}$, the time-dependent observables are then taken as
\begin{equation}
\A(\bm{x},t) = \A(\bm{\phi}(\bm{x},t))\label{Aeq}
\end{equation}
so that $\A(\bm{x},t)$ gives the value of the observables at time $t$ if the system was initially in state $\bm{x}$. The resulting time evolution of the observables can again be described by a differential equation
	\begin{equation}
    \frac{\partial}{\partial t}\A=%&
%	\sum_{i}\frac{\partial\phi_i}{\partial t}(\bm{x},t)\frac{\partial\A}{\partial\phi_i}(\bm{\phi}(\bm{x},t))\label{Ltime}\\
%	=&\sum_{ij}R_j(\bm{x})\frac{\partial\phi_i}{\partial x_j}(\bm{x},t)\frac{\partial\A}{\partial\phi_i}(\bm{\phi}(\bm{x},t))
%\nonumber
%\label{Ltime2}\\
%   =&\sum_jR_j(\bm{x})\frac{\partial}{\partial x_j}\A(\bm{\phi}(\bm{x},t))
%\nonumber\\=&
%= 
    L\A(\bm{x}),\hspace{0.8cm}
    L=\sum_iR_i(\bm{x})\frac{\partial}{\partial x_i} \label{Ltime}
	\end{equation}
with the Liouvillian $L$, a linear differential operator.
The general setup above requires us to track the full $\bm{x}$-dependence of the chosen observables $\bm{A}(\bm{x},t)$. To achieve a reduction in dimensionality, Chorin \cite{Chorin2000} assumes that the $\bm{x}$ are determined by $\bm{A}$ at least statistically, i.e.\ have some probability distribution that (only) depends on the current $\bm{A}$. Averages (expectations) over this distribution are written as $E[\cdot|\A]$, and the average evolution of $\bm{A}$ is governed by $\bm{v}(\A)= E[L\A(\cdot)|\A]$. Chorin \cite{Chorin2000,Chorin2002} showed that the corrections to this in the actual time evolution take the form of a memory term and a so-called random force $\bm{r}$, giving the 
%Putting this together we can obtain a 
general form for the time evolution of $\A$ as:
	\begin{equation}
	\frac{d}{d t}\A=\bm{v}(\A)+\int_0^tdt'\ \bm{M}(\A(t'),t-t')+\bm{r}\label{genEq}
	\end{equation}	
%We use $\bm{r}$ to represent the random force of the system. 
%We define $E[\cdot|\A]$ as the expectation of a function given a value of $\A$.
%So and describes the dynamics of $\bm{A}$ with itself and
The memory function 
$\bm{M}(\A(t'),t-t')$ depends on time difference $\tau=t-t'$ and -- nonlinearly -- on the past observable value $\A(t')$. Its evolution with $\tau$ is governed by the {\em deviations} of the drift from $\bm{v}(\A)$; this evolution reads for a general observable $g(\bm{x},\tau)$
\begin{equation}
%	\F(\bm{x})&=(L\A)(\bm{x})-E[L\A(\cdot)|\A(\bm{x})]\\
	\dfrac{\partial}{\partial\tau}g(\bm{x},\tau)=Lg(\bm{x},\tau)-E[(Lg)(\cdot,\tau)|\A(\bm{x})]\label{memTime}
\end{equation}
The memory function is obtained from the observable that measures exactly such fluctuations in the drift of the observables $\A$,
	\begin{equation}
	\bm{F}(\bm{x})=(L\A)(\bm{x})-E[L\A(\cdot)|\A(\bm{x})]
	\label{F_def}
	%\\\dfrac{\partial\F(\bm{x},\tau)}{\partial\tau}&=L\F(\bm{x},\tau)-E[(L\F)(\cdot,\tau)|\A(\bm{x})]\label{memTime}
	\end{equation}
From $\bm{F}(\bm{x})$ we define an $\bm{F}(\bm{x},\tau)$ 	 (\ref{memTime}) from the initial condition $\bm{F}(\bm{x},0)=\bm{F}(\bm{x})$, and the memory function is then given explicitly as $\bm{M}(\A,\tau) = E[L\bm{F}(\cdot,\tau)|\A]$. The random force itself is
$\bm{r}(\bm{x},t)=\bm{F}(\bm{x}.t)$ and has a vanishing average at all times, $E[\bm{r}(\cdot,t)|\A]=0$ \cite{Chorin2000}. In our context this term represents effects that come from the bulk starting \textit{away} from QSS. When the bulk starts in QSS, which is what we assume in the following, it vanishes and so can be discarded. While in \cite{Chorin2000} steady state dynamics are discussed, this is not required for the above formalism to be applicable.

% describes the dynamics of memory effects for $\A$. We use $\tau$ to represent time difference, which is present in the memory functions, as they integrate over past and current events. We then define in a similar way to \cite{Chorin2000}:

\subsection*{Subnetwork dynamics}

%General nonlinear equations

With the random force discarded as above, (\ref{genEq}) is a closed equation for the time evolution of the observables $\A$ and so achieves the desired dimensionality reduction. However, the memory function cannot in general be calculated in any closed form. We now show that this {\em can} be done, within a systematic approximation, for subnetwork dynamics. By this we mean that we consider as the observables $\A=\xs\label{subnetwork}$ a subset of $\bm{x}$, e.g.\ the concentrations of molecular species in a subnetwork of a large gene regulatory network. We denote the degrees of freedom in the rest of the network, the bulk, by $\xb{}$ and write out the components of the general time evolution (\ref{reaqsEx}) as
% as elements of $\bm{x}$ and term them as $\textrm{s}$ and $\textrm{b}$ respectively.
    \begin{equation}
	\frac{d\xs}{dt} = \rc{s}{(\xs,\xb{})},\hspace{0.8cm}
	\frac{d\xb{}}{dt} = \rc{b}{(\xs,\xb{})}\label{dxs_dt}
	\end{equation}
The Liouvillian then splits accordingly into
	\begin{equation}
	L=\sum_sR_s(\xs,\xb{})\frac{\partial}{\partial x_s}+\sum_bR_b(\xs,\xb{})\frac{\partial}{\partial x_b}
	\label{L_split}
	\end{equation}
	with the sums running over subnetwork and bulk species, respectively. Here and below subscripts always indicate individual species, while vectors with `s' and `b' superscripts collect all subnetwork and bulk quantities, respectively.
With (\ref{L_split}) the generic observable time evolution (\ref{Ltime}) $(\partial/\partial t)\xs = L\xs$ reduces to~(\ref{dxs_dt}) as it should. 
We now need to choose how to define the expectation $E[\cdot|\xs]$. We do this so that without the memory kernel, the reduced equation~(\ref{genEq}) corresponds to the simplification where the bulk dynamics equilibrates rapidly to any prevailing subnetwork state $\xs$, reaching a quasi steady state (QSS) value $\xb{*}$ defined by $d\xb{}/dt=0$ or
% at a "fast" speed, such that subnetwork time evolution depends only on its own concentration. This is represented by representing bulk dynamics with $\xb{*}$ which is defined as $\xb{}$ with a $0$ time derivative:
   	\begin{equation}
	\rc{b}{(\xs,\xb{*}(\xs))}=0\label{fastB}
	\end{equation}
As motivated in the introduction we will assume that this condition determines a unique bulk QSS $\xb{*}(\xs)$ for any $\xs$.  The expectation required to construct the reduced equation~(\ref{genEq}) is taken accordingly as 
\begin{equation}
E[g(\cdot) | \xs] = g(\xs,\xb{*}(\xs))\label{exp_def}
\end{equation}
i.e.\ by taking $\xs$ as prescribed and inserting for $\xb{}$ its QSS value.
The average drift $\bm{v}(\A) = E[L\A(\cdot)|\A]=E[\rc{s}(\xs,\xb{})|\xs]$ now evaluates directly from~(\ref{exp_def}) as
%  $E[L\A|\xs]$ will be of the form:
	\begin{equation}
\bm{v}(\xs) = 
%	E[L\A|\xs]=&E[\bm{R}(\xs,\xb{})|\xs]=
\rc{s}{(\xs,\xb{*}(\xs))}\label{subTime}
	\end{equation}
This is the QSS or `fast bulk' approximation to the subnetwork dynamics. Our main interest in the following lies in understanding the memory effects that account for the fact that the bulk is not in general fast, but evolves on a timescale comparable to that of the subnetwork.
To determine the resulting memory function, we start from the definition of $\bm{F}(\bm{x})$, which from (\ref{F_def}) has components
\begin{equation}
F_s(\xs,\xb{}) 
% &= R_s(\xs,\xb{}) - E[R_s(\cdot)|\xs]\\
=
R_s(\xs,\xb{}) - R_s(\xs,\xb{*}(\xs))
\label{F_time_zero}
\end{equation}
The main challenge is now to calculate the evolution of this observable in time according to~(\ref{memTime}). 
This is not feasible in general but we can develop a systematic approximation by {\em linearising} in deviations of the bulk degrees of freedom from the QSS, which we write as
% ie. the spirit of corrections to the This is generally a difficult problem to solve, we perform an approximation in order solve it. Our key approximation is that we will linearize $\F$ around the fast bulk steady state $\xb{*}$:
% definition of fast bulk species:
	\begin{equation}
	F_s(\xs,\xb{},\tau)\approx \sum_b(x_b-x_b^*(\xs))\fbs{b}{\tau}\label{memExp}
	\end{equation}
The problem then reduces to finding the evolution of $f_{bs}(\xs,\tau)$ from the initial condition 
$f_{bs}(\xs,0) \equiv f_{bs}^0(\xs)$ obtained by linearising~(\ref{F_time_zero}):
%** second equation could be inline to save space **
\begin{equation}
f^0_{bs}(\xs) = \frac{\partial R_s}{\partial x_b}
\label{f0_def}
\end{equation}
where the derivatives here and below are evaluated at $(\xs,\xb{*}(\xs))$ unless otherwise specified.

\subsection*{Memory evolution over time}

To derive our memory function we insert (\ref{memExp}) into (\ref{memTime}). Consistently applying the linearisation as detailed in Supp.~\ref{supp:expansion}  yields the following equation for $f_{bs}$: \begin{equation}   
	\frac{\partial}{\partial \tau}f_{bs}=
	\sum_{b'}\lbb f_{b's}+\sum_{s'}R_{s'}\frac{\partial}{\partial x_{s'}}f_{bs}\label{fdto}
\end{equation}
with
\begin{equation}
    \lbb %=\frac{\partial R_{b'}}{\partial x_b}-\sum_{s'}\frac{\partial R_{s'}}{\partial x_b}\frac{\partial x_{b'}^*}{\partial x_{s'}}
    =J_{b'b}+\sum_{s'b''}(\bm{J}^{-1})_{b'b''}\frac{\partial R_{b''}}{\partial x_{s'}}\frac{\partial R_{s'}}{\partial x_{b}}\label{l_def}
    \end{equation}
where the Jacobian matrix $\bm{J}$ is defined as
	\begin{equation}
	J_{b''b'}=\frac{\partial R_{b''}}{\partial x_{b'}}
	\label{Jac2a}
	\end{equation}
The next step is to find a solution $\fbs{b}{\tau}$ for the partial differential equation (\ref{fdto}). This can be done using the method of characteristics as the equation is linear in $\fbs{b}{\tau}$ and only involves first derivatives, and gives the closed form solution (see Supp.~\ref{supp:f})
% The solution effectively propagates the trajectory of the system forward in time and can be written explicitly as
%and We use a characteristic equation that propagates forwards in time the trajectory of the system:
\begin{equation}
	\fbs{b}{\tau}=\sum_{b'}E_{bb'}(\tau) 
	f_{b's}^0(\bm{\phi}_v(\xs,\tau))
	\label{fbs_sol}
\end{equation}
Here the $E_{bb'}$ are elements of the time-ordered matrix exponential $\bm{E}(\tau) = \exp[{\int_0^\tau d\tau'\,\bm{l}(\bm{\phi}_v(\bm{x}^\textrm{s},\tau'))}]$, and the propagation in time is performed with the flow $\phi_v$ for the QSS drift $\bm{v}(\xs)$.   

\subsection*{Memory function}
We can now finally determine the memory function on subnetwork species $s$, which from the general framework set out above is
\begin{equation}
	M_s(\xs,\tau)=E[LF_s(\cdot,\tau)|\xs]\label{eq:mem}
\end{equation}
We insert the expansion (\ref{memExp}) here and obtain after some algebra (see Supp.~\ref{supp:expansion}) our main result, a simple expression for the memory function:
\begin{equation}
M_s(\xs,\tau)
%E[LF_s(\cdot,\tau)|\xs]
=
\sum_{b'}c_{b'}(\xs)f_{b's}(\xs,\tau)
\label{M_final}
\end{equation}
where we have denoted
\begin{equation}
c_{b'}(\xs) = \sum_{s'}
%\left(
\sum_{b''}(\bm{J}^{-1})_{b'b''}\frac{\partial R_{b''}}{\partial x_{s'}} R_{s'} \label{c_def}
%\right)
\end{equation}
%For the sake of practicality we define:
%\begin{align}
%	c_{b'}(\xs)=\sum_{s'}R_{s'}\left(\sum_{b''}(J^{-1})_{b'b''}\frac{\partial R_{b''}}{\partial x_{s'}}\right)
%\end{align}
These functions 
%Where $\bm{c}(\xs)$
can be thought of as prefactors to the memory term.
The general projected time evolution equation now takes the form
	\begin{equation}  
	\frac{d}{d t}x_s=v_s(\xs(t))+\int_0^tdt'\ M_s(\xs(t'),t-t')%+r_s(t)
%\label{eqAbs}
\label{subEq}
    \end{equation}  
%We can then write (\ref{eqAbs}) in the form:
%\begin{align}
%	\frac{\partial}{\partial t}x_s&=R_s(\xs,t)+\int_0^tdt'\sum_{b'}c_{b'}(\xs(t'))f_{b's}(\xs(t'),{t-t'})\label{subEq}
%\end{align}
The first term contains the QSS drift while the second one represents the memory correction to this, which is expressed in terms of the memory function (\ref{M_final}). Our derivation allows this memory to cover the behaviour around multiple fixed points of the system, due to its nonlinear dependence on $\xs$. The interpretation of our result~(\ref{M_final}) is that in a small time interval $dt'$, $x_b-x_b^*$ will change by $c(\xs(t'))\,dt'$. This deviation from the QSS is propagated by the exponential matrix and affects the drift $R_s$ at time $t$ as captured by  $f^0_{bs}$ in (\ref{fbs_sol}).
In Supp.~\ref{sSec:Gouasmi} we compare (\ref{M_final}) with the work of \cite{Gouasmi2017}, which instead of the QSS assumption takes $\xb{}=0$.
%the bulk being zero (.
This %assumption
is unsuitable for the multistable systems we are interested in but we show that the method can be adapted to project to bulk QSS values (Supp.~\ref{sSec:GQSS}). This leads
%find that when mapping such method to bulk QSS values we obtain 
to an expression similar to (\ref{M_final}), but crucially without the propagation in time from $t'$ to $t=t'+\tau$ % The difference between the two methods consists in that the method from Supp.~\ref{sSec:GQSS} does not propagate the deviations over $\tau$ (\ref{M_g})
(Fig.~\ref{sFig:bisMem}).

\subsection*{Self-consistent approximation}
Our linearisation approach (\ref{memExp}) implies that the memory term captures dynamical effects that are of first order in the deviations of the bulk network from its QSS. We will now develop an approximate self-consistent way of incorporating higher order corrections, which turns out also to simplify the numerical evaluation of the memory terms. Consider the factor $f^0_{bs}(\bm{\phi}_v(\xs,\tau))$ that from (\ref{fbs_sol},\ref{M_final}) appears in the memory function $M_s(\xs,\tau)$. In the actual memory integral this is evaluated for $\xs(t')$ and $\tau=t-t'$, i.e.\ as $f_{bs}^0(\bm{\phi}_v(\xs(t'),t-t'))$.   
As explained above, $\phi_v$ is the flow generated only by the QSS drift, i.e.\ without memory corrections. But the memory terms change the flow, so we can make the approach self-consistent by substituting for $\phi_v$ the actual time evolution {\em with} memory. This corresponds to replacing
\begin{equation}
\bm{\phi}_v(\xs(t'),t-t') \to \xs(t)
\end{equation}
as we are just propagating the subnetwork state from $\xs(t')$ by a time difference $t-t'$ to $\xs(t)$. Making this replacement also in the matrix exponential in (\ref{fbs_sol}) changes the memory term $\mathcal{M}_s(t) = \int_0^tdt'\ M_s(\xs(t'),t-t')$ into
\begin{align}
\tilde{\mathcal{M}}_s(t) =& \sum_{b''}\int_0^{t}dt'\ \sum_{b'} c_{b'}(\xs(t'))\left(e^{\int_{t'}^t dt''\,\bm{l}(\xs(t''))}\right)_{b'b''}\nonumber\\
&\times f_{b''s}^0(\xs(t))\label{M_selfconsistent}
\end{align}
The dependence on the subnetwork species $s$ on which the memory acts is contained only in the -- now $t'$-independent -- factor $f_{b'',s}^0(\xs(t))$. As shown in Supp.~\ref{supp:selfcon}, the memory integrals in the first line can then be calculated efficiently as solutions to differential equations, one for each bulk species $b''$. Conceptually, however, the self-consistent memory term is more complicated. In the original formulation (\ref{subEq}), the memory is a superposition of separate effects from all past times $t'$: the state $\xs(t')$ of the subnetwork affects the behaviour of the bulk and feeds back into the subnetwork at time $t$. In (\ref{M_selfconsistent}), the way this feedback acts is additionally modulated by the entire time evolution of the subnetwork between times $t'$ and $t$. In the applications considered below both approaches yield similar quantitative results, hence which one to choose depends on the aim: for numerical calculations of memory effects the self-consistent version is more efficient, whereas the memory functions themselves are easier to analyse in the original version because they depend -- in addition to time difference, which always features -- only on the subnetwork state at one time $t'$. Note that what we refer to as self-consistency is distinct from an approach widely used in application of projection methods to physical systems (see e.g.~\cite{Kawasaki1970,Gotze1991}), where equations of motion for correlation functions are considered and the relevant memory kernels are, via appropriate approximations, related back self-consistently to the correlation functions.

\subsection*{General memory properties}

Both the \zmv and \zms expressions for the memory term that we have derived are, as we have emphasized, nonlinear in $\xs$ and so not of the convolution form that would appear in linear ZM projection methods. This simpler form is recovered, however, in the dynamics near fixed points. To see this, we note that at a global fixed point $(\xs{}^*,\xb{*})$, where $R_s=0$,
 %the bulk is necessarily at its QSS, $\xb{}=\xb{*}(\xs)$, given our assumption that the bulk always has a unique steady state given the subnetwork. Thus
the last factor in the definition (\ref{c_def}) of $c_{b'}$ 
%, $R_s(\xs,\xb{*})$, equals $R_s(\xs,\xb)$ and so 
vanishes. As $c_b$ is a factor in both of our memory expressions, both receive zero contributions when $\xs(t')$ is at a fixed point. Assuming $R_s$ and $R_b$ are sufficiently smooth, we can therefore {\em linearize} the memory terms for dynamics near such a fixed point. In this linearization all other factors in the memory are evaluated at the fixed point, and from this one deduces (see Supp.~\ref{supp:linear}) that the linearized forms of the \zmv and \zms memory expressions are in fact identical. The memory terms then become time convolutions of $\xs(t')-\xs{}^*$ with a memory kernel that depends only on time; Supp.~\ref{supp:linear} provides a numerical example. It follows from the general arguments in Supp.~\ref{supp:exact} that these kernels describe {\em exactly} the dynamics of the full linearized system, and so the equations with memory will correctly predict e.g.\ the relaxation behaviour near a stable fixed point.
% An example application of the linearized memory analysis is shown in Supp.~\ref{supp:linear}.
%When the system is sufficiently near a fixed point, the dynamics will behave in a linear fashion. The memory functions can be justifiably linearised, resulting in the memory \zmv{} and \zms{} becoming equivalent convolution-type equations. Thus when near a fixed point the projected equations will predict a relaxation at long time scales.

\subsection*{Memory decomposition}

In spite of their nonlinearity it turns out to be possible to decompose our memory expressions into specific ``channels'', in order to analyse the contribution of interactions within a network. Building on the approach of \cite{Herrera-Delgado2018}, we take advantage of the two partial derivative expressions in (\ref{f0_def}, \ref{c_def}) to decompose the memory exactly into combinations of incoming and outgoing channels (Supp.~\ref{sSec:decomp}).
The analogous construction for the self-consistent approximation is set out in 
%follows a similar construction and thus we can also perform a channel decomposition there (
Supp.~\ref{sSec:selfConDec}.

\section*{Applications \label{sec:appCh}}
To test the effectiveness of the method we examine systems that contain multiple steady states, oscillatory behaviours and complex transient dynamics. These are relevant in a wide range of physical and biological contexts.

\subsection*{Multistability}
We first examine a series of multistable systems defined by mutually repressive Hill functions: %** this is "or"-repression, i.e. production of $x_i$ is inhibited if either of the neighbours has large $x_j$. Comment on this? Is it a standard choice? Adding separate Hill functions for each $j$ would give an "and" repression **
\begin{equation}
\frac{d}{dt}x_j=\frac{a}{1+\sum_{i\neq j}x_i^n}-x_j \label{eq:multistable}
\end{equation}
The above equation constitutes an ``or'' logic because of the sum of the terms in the denominator, where even if only one repressor has a high concentration, the production rate will become very low. Such interactions lead to multistability in a wide variety of developmental systems \cite{Angeli2004}.

\begin{figure}[h]
\centering
% \includegraphics[width=0.35\linewidth]{bistable_net.pdf}\\
% \includegraphics[width=0.58\linewidth]{bistable_memTD.pdf}
% %\includegraphics[width=0.22\textwidth]{USS_a=4_n=2_halfSize_2.pdf}
% \includegraphics[width=0.4\linewidth]{USS_a=4_n=2_halfSize.pdf}
\includegraphics[scale=0.5]{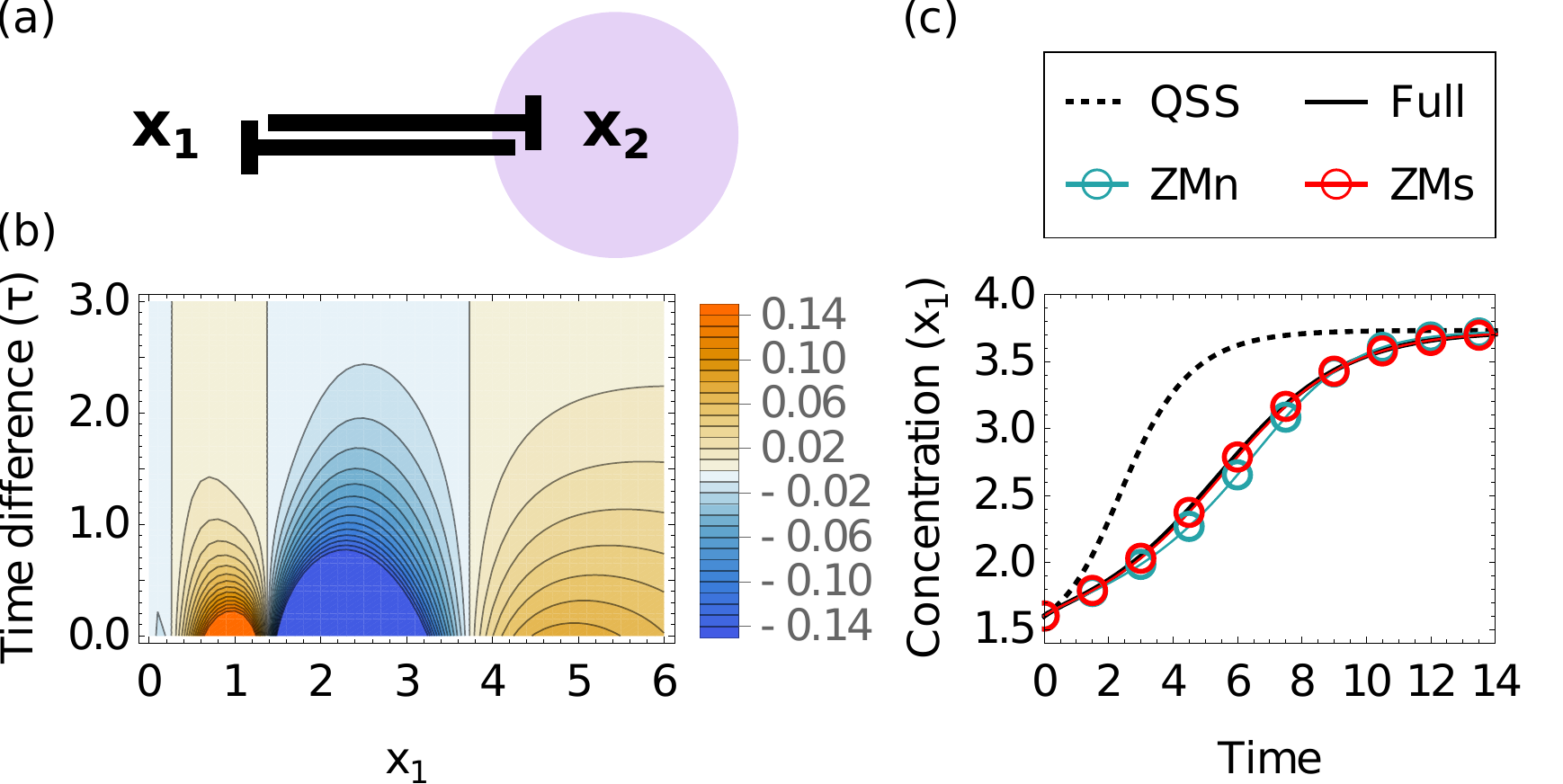}
\caption{(a) Network illustrating a bistable switch defined using cross-repressive Hill functions (\ref{eq:multistable}) with $a=4$, $n=2$ (steady states are $(x_1,x_2)=(6,0.028)$ and $(0.028, 6)$ with an unstable fixed point at $(1.46,1.46)$). For this and all other network illustrations, blunt arrows indicate repression; purple shading identifies the species placed in the bulk. (b) Memory function for the bistable switch shown as a function of (past) concentration $x_1$ of species 1 ($x$-axis) and
% axis indicates the concentration of $x_1$ and $y$ axis indicates how the memory affects $x_1$ over time after the corresponding $x_1$ concentration occurred. 
time difference $\tau$ ($y$-axis).
Memory function values range from negative to positive as indicated by the scale bar on the right, and are capped by blue and orange outside the scale bar range. (c) Time course of the system demonstrates the capacity of both the nonlinear (ZMn; cyan) and self-consistent (ZMs; red) projections to capture the timescale and shape of transients of the full model (solid line) in reaching a stable fixed point. The QSS approximation (dashed line) significantly underestimates the length of the transient, showing that the ZM projections successfully correct for $x_2$ not being at QSS. \label{fig:memShape} \label{fig:bistableNet} \label{fig:multSS}}
\end{figure}

We test the method on the simplest case with two nodes $\{x_1,x_2\}$, which leads to a system that cannot produce oscillations \cite{Page2018}. We place $x_2$ in the bulk and calculate the memory function for the single remaining subnetwork species $x_1$ (Fig.~\ref{fig:bistableNet}a). This depends on the past concentration $x_1(t')$ and the time difference $\tau=t-t'$ (Fig.~\ref{fig:memShape}b). We observe that the memory becomes zero at each fixed point as expected from \cite{Herrera-Delgado2018}, where the memory was obtained as an expansion (to quadratic order) in devations of $\xs$ from a fixed point. To leading order the memory grows linearly with this deviation, and in line with this we see it changing sign at every fixed point. The sign of the memory in all cases is opposite to that of the drift, so the memory delays the relaxation time to the corresponding steady state. This makes intuitive sense as in the original system, the bulk species' state reacts relatively slowly to subnetwork changes, rather than infinitely fast as the QSS approximation assumes.

To test the accuracy of our method in capturing the transient temporal dynamics we set the initial condition of $x_1$ to be close to the unstable fixed point of the QSS dynamics; here we are furthest from the stable fixed points and so can test the limits of the method. For the evaluation we used as a baseline the full dynamics of the original system, setting $x_2$ at time zero to its QSS value with respect to the value of $x_1$. We compare this to the subnetwork dynamics predicted by the simple QSS dynamics, and by our approach, which includes memory corrections. For this example we evaluate both the nonlinear memory description (\ref{subEq}), which we label \zmv (Zwanzig-Mori nonlinear), and the self-consistent memory (\ref{M_selfconsistent}),
denoted \zms (Zwanzig-Mori self-consistent) below. We find that both replicate the behaviour of the original system well, independently of whether the initial condition eventually leads to the low- or high-$x_1$ fixed point. The QSS approximation, on the other hand, reaches the steady state unrealistically fast (Fig.~\ref{fig:multSS}c). Given that the \zms is substantially easier to implement for time course prediction (Supp.~\ref{supp:selfcon}) we concentrate on this approach below. Further justification for this comes from the fact that the self-consistent memory description is {\em exact} when $R^{\rm s}$ and $R^{\rm b}$ depend at most linearly on the bulk species, as we show in Supp.~\ref{supp:exact}. This exactness is not a trivial consequence of the fact that our approach is a linearisation in $\xb{}-\xb{*}$, as it would otherwise hold also in the ZMn version. Biochemical systems with linear $\xb{}$-dependences usually involve mass-action reactions and can produce bistable systems or oscillations \cite{Wilhelm2009,Kondepudi1998}, which we can then reproduce exactly with the \zms projection (see Figs.~\ref{sFig:exactOsc},~\ref{sFig:exactBis}). A well-known example from physics is the Caldeira-Leggett model of a heat bath. This has a bulk composed of harmonic oscillators\cite{Caldeira1983}, and the resulting memory term was obtained exactly by Zwanzig \cite{Zwanzig1973,Zwanzig2001}. Since in that system the bulk degrees of freedom appear linearly, our approach reproduces this solution but is more general in that it does not, for example, rely on the dynamics to be derived from a Hamiltonian.

\begin{figure}[h]
\centering
\includegraphics[scale=0.5]{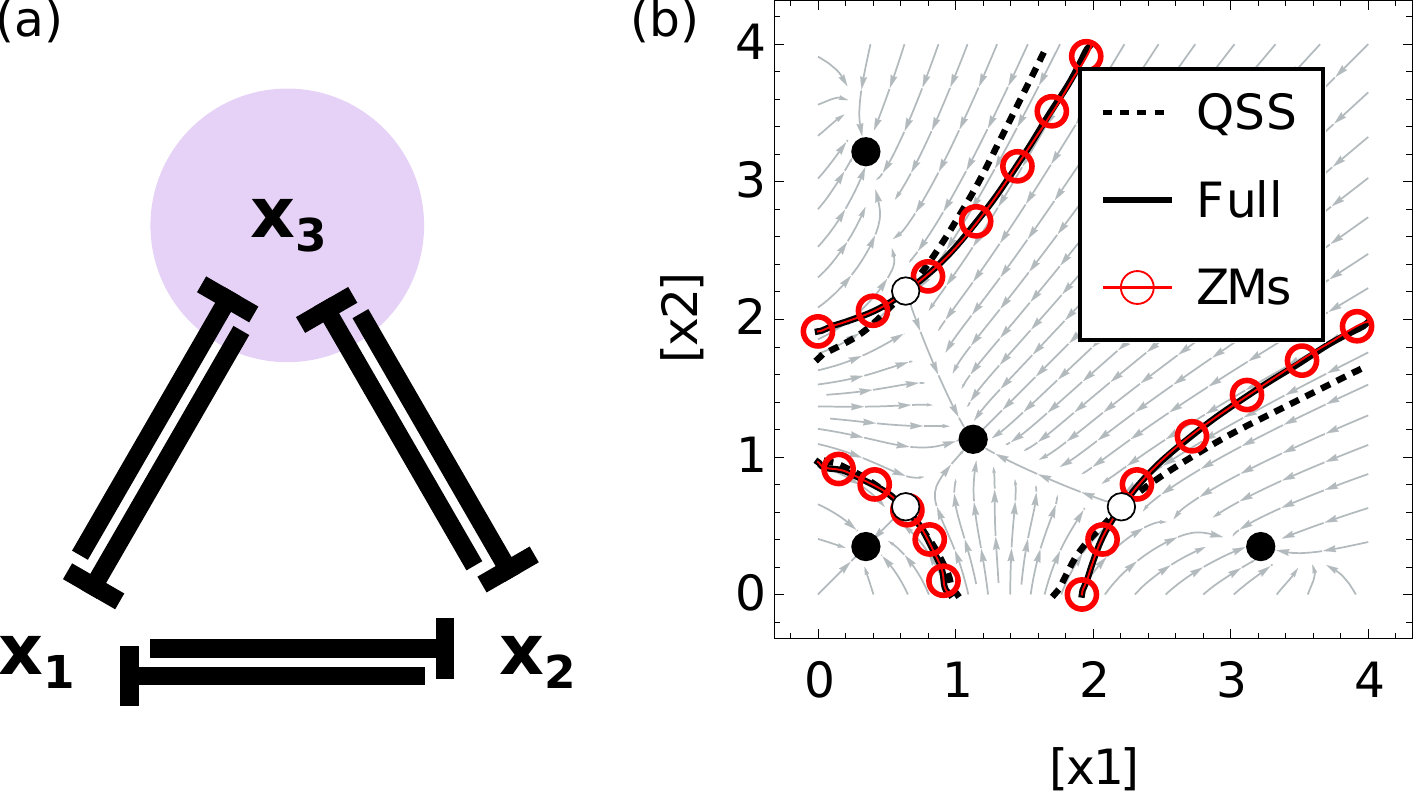}
\caption{(a) Network illustrating the cross-repressive tetrastable system (\ref{eq:multistable}); purple shading indicates the species placed in the bulk. Parameters are $a=4$, $n=2$.
% (steady states are ${x_1=3.219, x_2=0.348, x_3=0.348}$, $\{x_1=0.348, x_2=3.219, x_3=0.348\}$, $\{x_1=3.219, x_2=0.348, x_3=0.348\}$ and $\{x_1=1.128, x_2=1.128, x_3=1.128\}$ with unstable fixed points $\{x_1=2.206, x_2=0.638, x_3=0.638\}$, $\{x_1=0.638, x_2=2.206, x_3=0.638\}$ and $\{x_1=0.638, x_2=0.638, x_3=2.206\}$). 
(b) Phase portrait indicating the basins of attraction of the four stable fixed points (black circles) and the unstable fixed points (white circles). The separatrices bounding each basin are shown for the full dynamics (solid lines), QSS approximation (dotted lines) and the subnetwork equations with memory (ZMs, red circles); stream plots are shown for the QSS approximation. The QSS approach shows a clear difference to the original system, while the boundaries set by \zms and the full system are almost indistinguishable. \label{fig:triMultSS} \label{fig:tetrastableNet}}
\end{figure}

We next tested the approach on a tetrastable system defined in the same way as (\ref{eq:multistable}) with variables $\{x_1,x_2,x_3\}$. We consider the subnetwork containing $x_1$ and $x_2$, which will allow us to investigate the effect of the memory effects on the shapes of the basins of attraction of the different (stable) fixed points (Fig.~\ref{fig:tetrastableNet}a). For the parameter values we use, there are four such fixed points for the full network: three where only one species has high concentration and the other two low, and one where all concentrations are equal (Fig.~\ref{fig:tetrastableNet}b).
The boundaries of the basins of attraction can be read as the points where, depending on its initial condition, the system chooses a different basin of attraction. In a biological setting these choices could represent cell fate decisions, where a cell decides to adopt a specific specialised function and become a particular cell type (for a review see: \cite{Enver2009}). We find that the QSS system fails to replicate the decision process of the original system, whereas the \zms accurately identifies both the eventual steady state (Fig.~\ref{fig:triMultSS}b) \textit{and} the timing to get to this state (data not shown).

\subsection*{Oscillations}
\begin{figure}[h]
\centering
\includegraphics[scale=0.5]{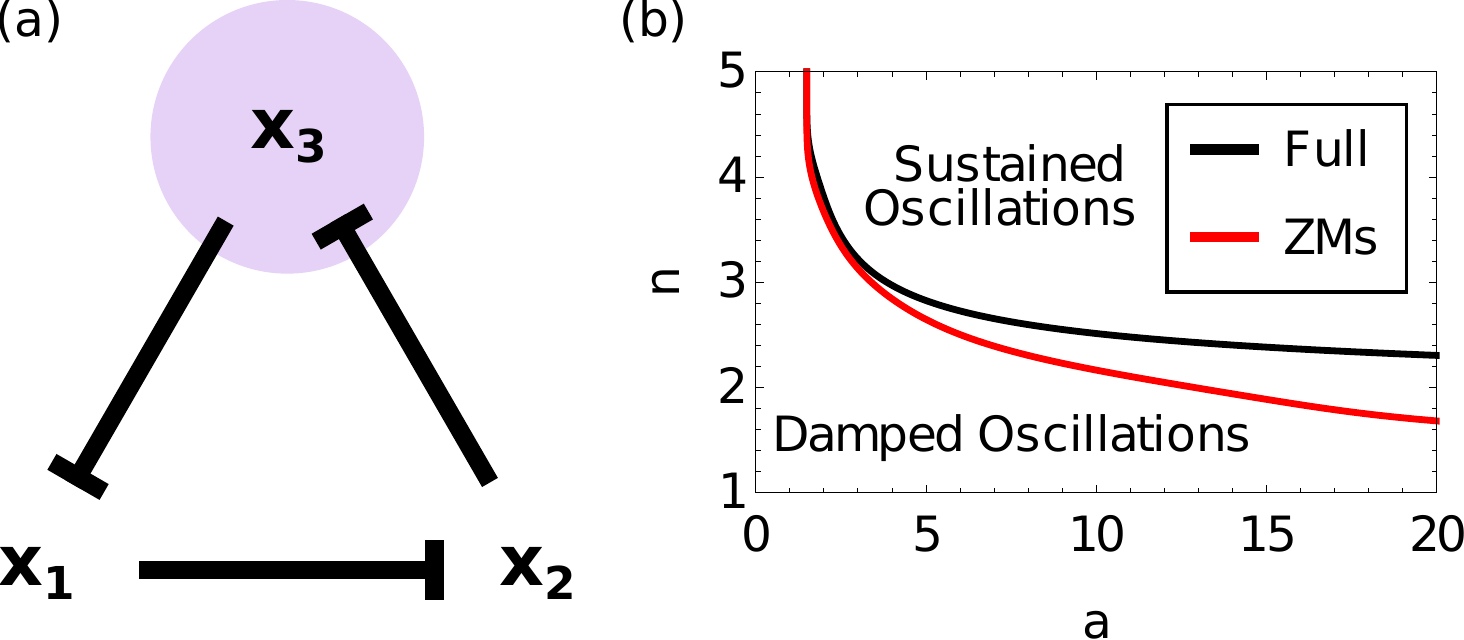}
\caption{(a) Network illustrating the repressilator system (\ref{eq:repressilator}); purple shading indicates the species placed in the bulk. (b) Bifurcation diagram of the repressilator for system parameters $a$ and $n$. The lines represent super-critical Hopf bifurcations. The QSS system can only produce damped oscillations, so has no bifurcation at all. The \zms system (red line) shows a good qualitative match to the shape and position of the bifurcation of the full system (solid black line). \label{fig:repressilatorNet} \label{fig:bifOsci}}
\end{figure}

We further explore the ability of the subnetwork equations with memory to reproduce oscillations arising from a uni-directional repressive network. We use the repressilator system \cite{Elowitz2000}, which robustly generates oscillations due to delays arising from an odd number of nodes \cite{Strelkowa2010,Page2018}. It has concentration variables $\{x_1,x_2,x_3\}$ and repressive interactions as shown in Fig.~\ref{fig:repressilatorNet}a and represented mathematically by
\begin{equation}
\frac{d}{dt}x_j=\frac{a}{1+x_{j-1}^n}-x_j
\label{eq:repressilator}
\end{equation}
where $x_0 \equiv x_3$.
We first compare the bifurcation diagram that results from varying both system parameters $a$ and $n$, in order to see whether the subnetwork equations with memory can replicate the 2D Hopf bifurcation of the original system, from damped to sustained oscillations (Fig.~\ref{fig:bifOsci}b).

\begin{figure}[h]
\centering
\includegraphics[scale=0.5]{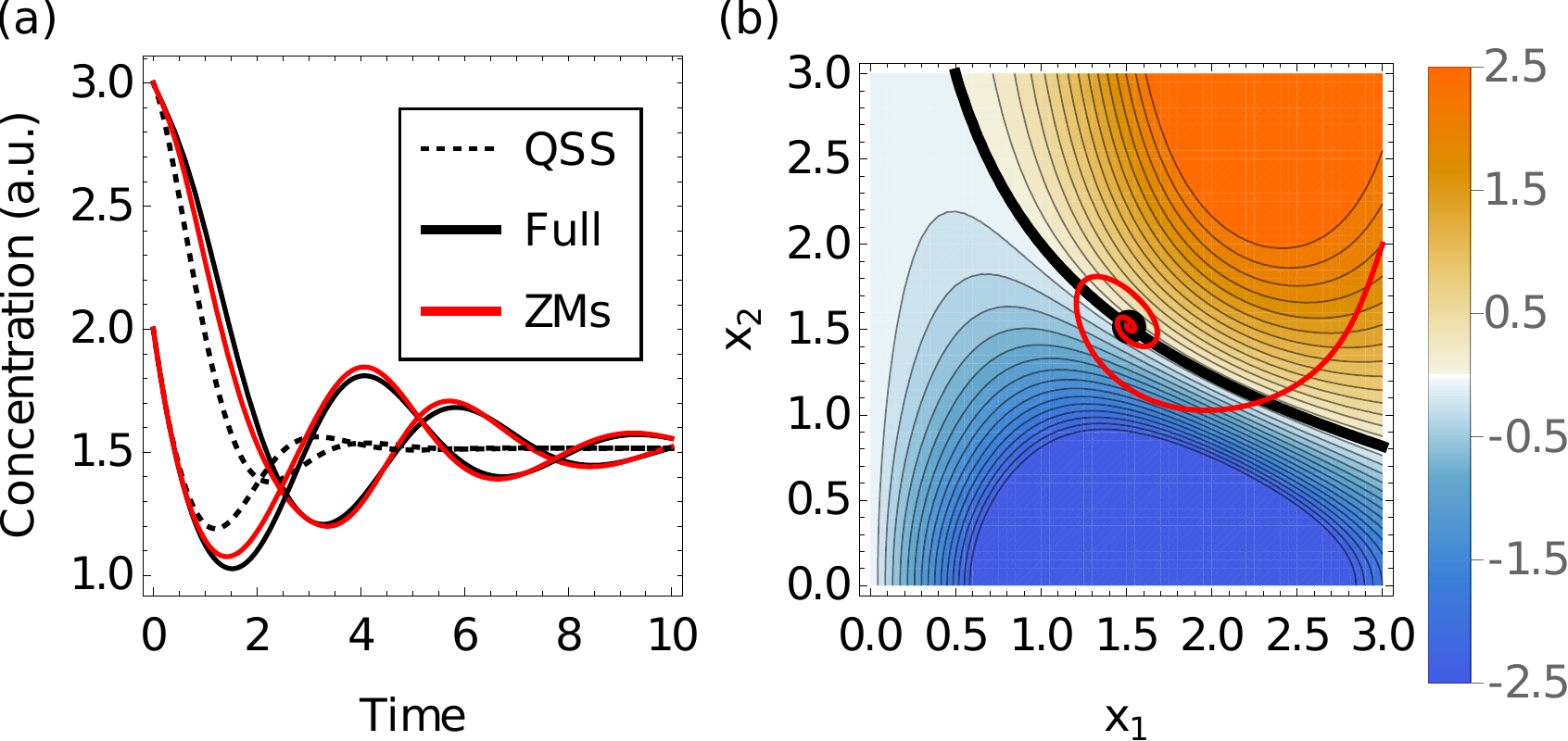}
\caption{ (a) Damped oscillations generated by the repressilator \cite{Elowitz2000} can be accurately reproduced by the projection approach (red line) while the QSS (dotted line) fails to replicate both the timing and the amplitude of the oscillations of the full system (solid black line). (b) Colour map of memory amplitude (memory at $\tau=0$) as a function of ($x_1$, $x_2$).
The memory amplitude changes from negative (blue) to positive (orange) across the thick black line. % indicates where the memory value is zero.
Red curve: Parametric plot of the \zms time course from the left;
% is represented by the red curve,
the memory repeatedly changes from negative to positive to drive the correct oscillations.
\label{fig:dampMemA}
\label{fig:dampOsci}}
\end{figure}

%While the QSS system is unable to generate sustained oscillations at all, 
In contrast to the QSS approximation we find that the projection technique correctly replicates the existence of sustained oscillations and predicts a qualitatively correct bifurcation diagram. The period and amplitude of sustained oscilations in the relevant parameter regime is less well replicated (not shown). For damped oscilllations the subnetwork equations with memory work accurately in predicting the full temporal dynamics (Fig.~\ref{fig:dampOsci}a). By contrast, the QSS approximation displays almost no oscillatory behaviour and none of a sustained nature, highlighting the importance of memory effects for oscillatory transients.

In order to understand in more detail how memory generates oscillations, we analyse the corresponding memory functions. We first plot the memory amplitude, i.e.\ the value $M(\xs,0)$ for memory from the immediate past ($t'=t$) across the configuration space of our subnetwork (Fig.~\ref{fig:dampOsci}b) and observe two distinct regions with positive and negative memory amplitude, separated by a line where this amplitude vanishes (black). Plotting the time course from Fig.~\ref{fig:dampOsci}a in the same representation we observe that it crosses the black line many times. The corresponding changes in the sign of the memory amplitude %changes sign approximately at the same time and so
are what drives the oscillations seen in Fig.~\ref{fig:dampMemA}a.

\section*{Neural tube network (transients and multistability) \label{sec:bioCh}}
\begin{figure}[h]
\centering
\includegraphics[scale=0.5]{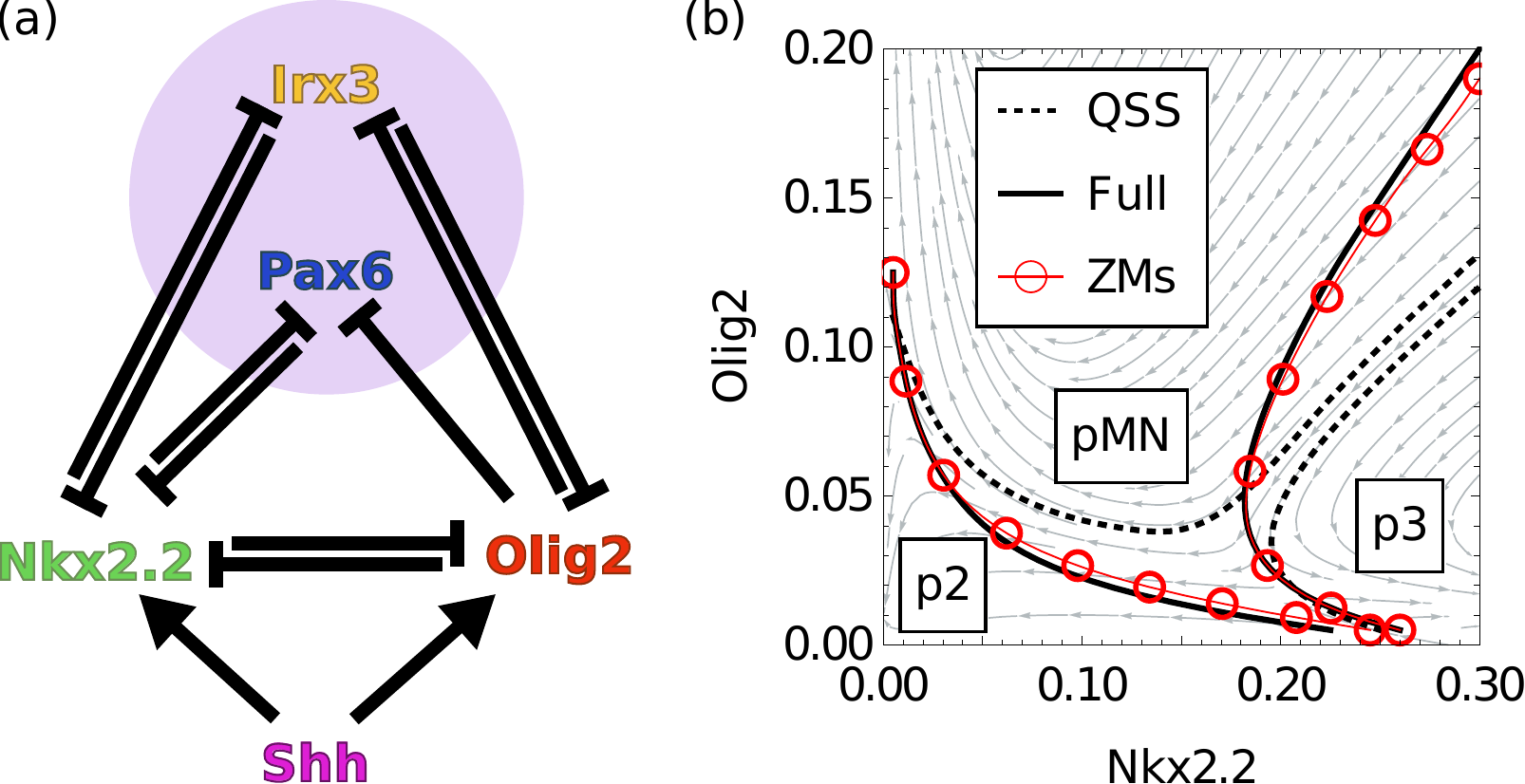}
\caption{(a) Neural tube network \cite{Cohen2014} defined by cross-repressive interactions between four transcription factors and activation by Sonic Hedgehog signalling (Shh), dependent on neural tube position; purple shading indicates choice of bulk. (b) Fate decision diagram for a neural tube position with three attractors (at position $p=0.65$); solid lines indicate boundaries of basins of attraction that biologically separate different fate choices. Possible steady states are p3 (high Nkx2.2, right), pMN (high Olig2, top) and p2 (low Nkx2.2 \& Olig2, and high Irx3 \& Pax6, bottom left). Dashed lines indicate basin boundaries for the QSS approximation and red dots the basin boundaries for the \zms projection; a stream plot is shown for the QSS system. The \zms system very accurately reproduces the boundaries of the full system. \label{fig:PONI_tristable}}
\end{figure}

Finally we apply the projection approach to a biologically relevant system with several bifurcations and non-trivial dynamical properties, specifically the neural tube gene regulation network described in \cite{Cohen2014}. The network is sketched in Fig.~\ref{fig:PONI_tristable}a and its specific form and parametrisation is given in Supp.~\ref{sSec:NT}. In each cell of the neural tube the network responds to the concentration of a extracellular signalling molecule called Sonic Hedgehog. 
This molecule forms a gradient in the developing vertebrate neural tube which is interpreted by cells in the tissue through the actions of a GRN (for a recent review see \cite{Sagner2019}).
As a result of the concentration changing systematically along the neural tube, we are in fact dealing here with a family of networks that vary with neural tube position, parametrised below in terms of $p$ running from 0 to 1 (where $p$ also represents the Shh input).
% where the $p$ represents the Shh input.
The network contains four molecular species (called transcription factors in the gene regulation context), two of which generate a bistable switch by mutual cross-repression as in our first example above. Similarly to how we proceeded in \cite{Herrera-Delgado2018}, we therefore place these two molecular species (Nkx2.2 and Olig2) in the subnetwork, with the bulk provided by the two other species (Irx3 and Pax6).

We test the method at the position along the neural tube where the model has the most complexity, a region of tristability ($p=0.65$), and compare with the original system and the QSS approximation (Fig.~\ref{fig:PONI_tristable}b). We find that as for the tetrastable example above (Fig.~\ref{fig:triMultSS}), the projection accurately replicates the choice of steady state, in contrast to the QSS method (Fig.~\ref{fig:PONI_tristable}b). With the memory included, three basins of attraction are predicted, of which one - labelled pMN - separates the other two (p2 and p3) so that no direct transitions from p2 to p3 can occur. The QSS approximation is not just quantitatively inaccurate but loses this biologically important qualitative feature, predicting instead that the p2 and p3 basins border each other. At other neural tube positions we also consistently find a good match between the original system and the \zms projection approach (not shown).

We next analyse the temporal evolution of the systems at various neural tube positions, using the experimentally determined initial condition for Nkx2.2 and Olig2 (vanishing concentration); we again compare the \zms description with the original system and the QSS reduction. The \zms predictions show a good fit with the original system at all positions (Fig.~\ref{fig:PONI_transient} displays results for a position with a strongly non-monotonic transient, $p=0.1$), demonstrating that the memory functions are capable of accurately capturing not just final cell fate decisions but also the timing of such decisions. This temporal aspect %of patterning 
is important for correct patterning as explored in \cite{Exelby2019}.

\subsection*{Decomposing nonlinear memory functions}

\begin{figure}[h]
\centering
\includegraphics[scale=0.5]{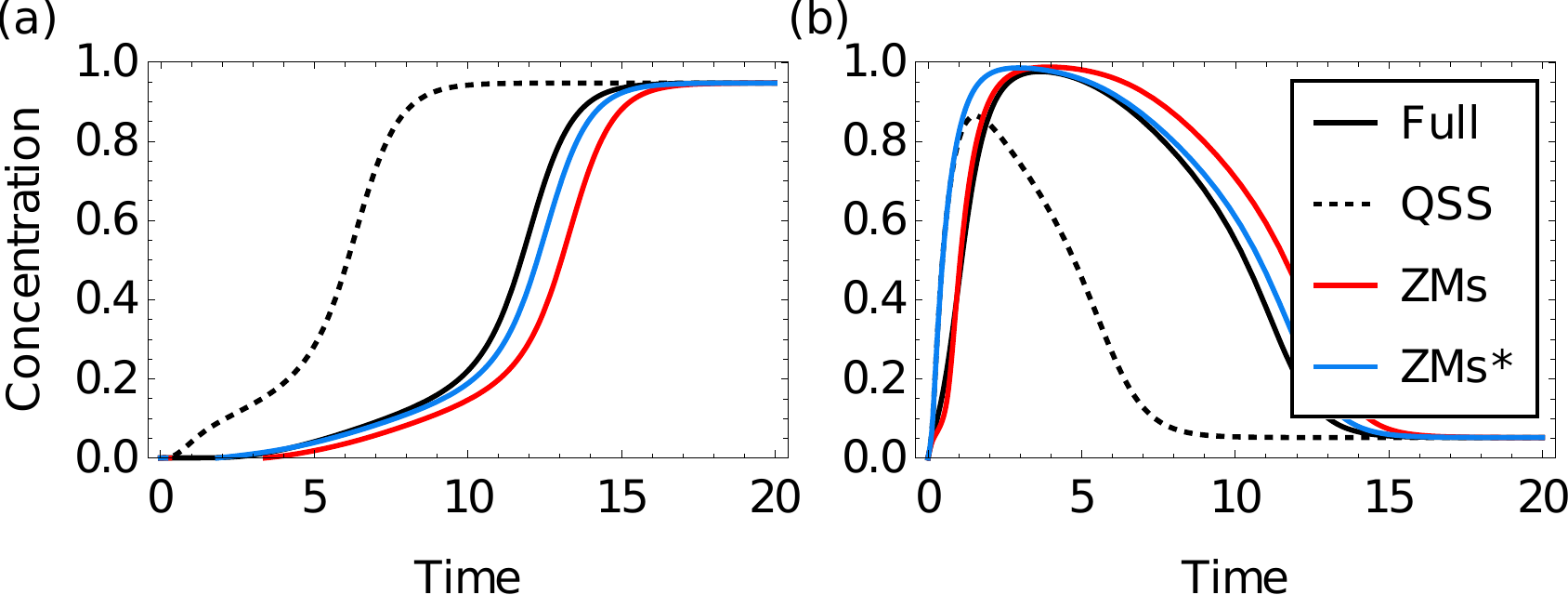}
\caption{Time courses of concentrations of Nkx2.2 (a) and Olig2 (b) in p3 domain ($p=0.1$). A transient expression of Olig2 leading to a delay in Nkx2.2 expression is observed \textit{in vivo}. The full system (solid line) and \zms projections (with Nkx2.2 and Olig2 chosen as subnetwork, see Fig.~\ref{fig:PONI_tristable}a, red and blue) qualitatively reproduce this behaviour. In contrast, the QSS approximation (dashed line) is unable to capture the long Olig2 transient. ZMs$*$ represents the removal of all memory functions except those specified in Fig.~\ref{fig:mMem}a, which suffice to capture the observed transient. \label{fig:PONI_transient}}
\end{figure}

\begin{figure}[h]
\centering
\includegraphics[scale=0.5]{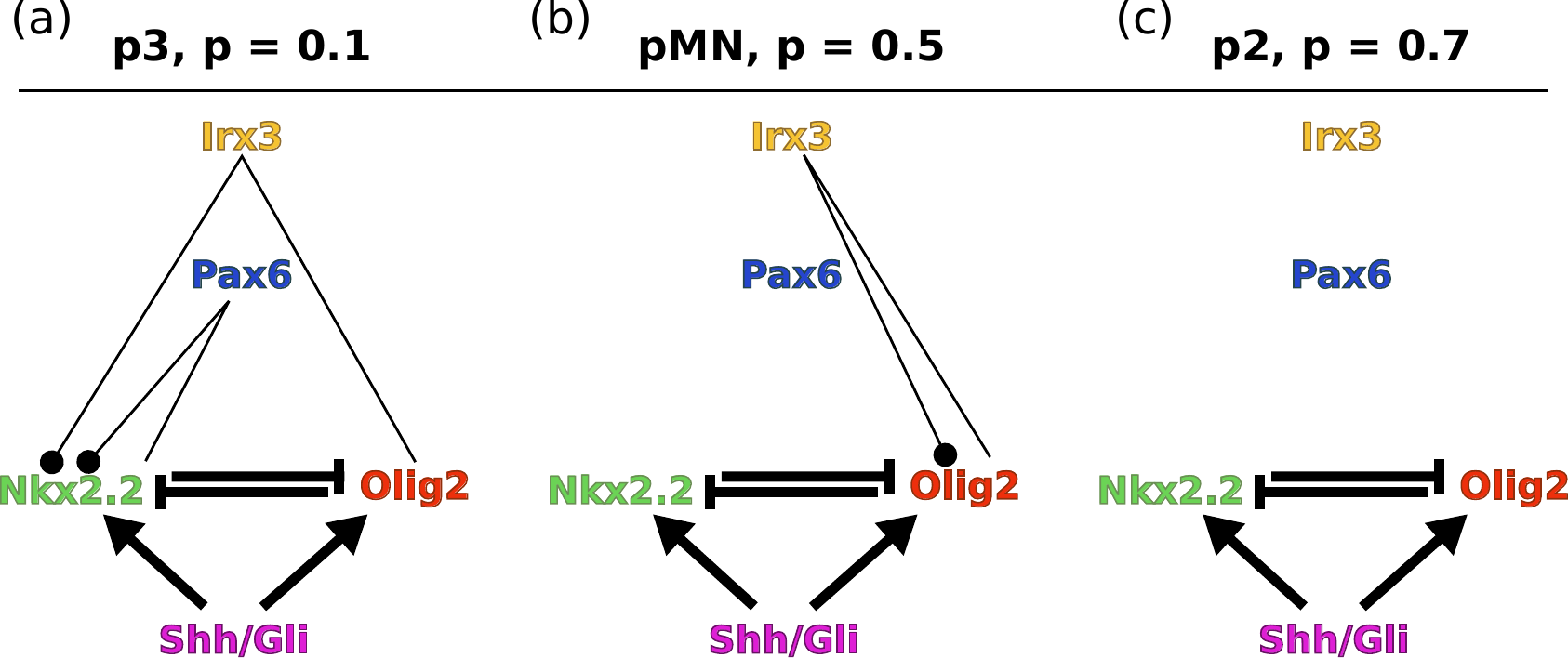}
\caption{Diagrams indicating the channels that have the largest memory contributions to the observed dynamics at three distinct neural tube positions. Black dot indicates the species receiving the memory contribution; the other end of each line is the species ``sending'' the memory. Contributions change according to the final steady state: p3 (a), pMN (b) and p2 (c). \label{fig:mMem}}
\end{figure}

In order to understand how memory functions affect the patterning dynamics, we set out to understand their structure. We perform our self-consistent memory decomposition approach (Supp.~\ref{sSec:selfConDec}) and analyse the results to identify the memory channels with the most impact on the time courses based on their contribution along the trajectory (Fig.~\ref{sFig:decompP3}). Performing this analysis for the different progenitor domains predicts the most important regulatory interactions contributing to the memory effect at each neural tube position (Fig.~\ref{fig:mMem}). This indicates marked differences in the most significant memory channels at different neural tube positions (see Fig.~\ref{sFig:decompP3} for an illustration of the decomposition in the p3 domain at $p=0.1$).

To test the validity of our results we remove all channels identified as unimportant, setting them to QSS, thereby keeping only the channels shown in Fig.~\ref{fig:mMem}a. Simulating the dynamics with only these memory functions results in dynamics that match closely those of the full simulation (Fig.~\ref{fig:PONI_transient}), confirming the prediction that these channels dominate the memory effects.

We next investigated the experimentally validated transient in gene expression in the p3 domain (monostable with high Nkx2.2 in the steady state, $p=0.1$, Fig.~\ref{fig:PONI_transient}). Nkx2.2 induction is delayed in neural progenitor cells compared to Olig2 \cite{Dessaud2007} and our analysis of the memory function provides insight into how this is achieved. The active memory channels ensure that Nkx2.2 is kept close to zero while Olig2 rises (Fig.~\ref{fig:PONI_transient}). The dominant memory channels shown in Fig.~\ref{fig:mMem}a indicate that a different bulk species captures the history of each subnetwork species: Pax6 transmits the memory of Nkx2.2 and Irx3 the one of Olig2. The effect of these these bulk species is thus to delay Nkx2.2 expression based on the past expression of Nkx2.2 and Olig2.

\begin{figure}[h]
\centering
\includegraphics[scale=0.5]{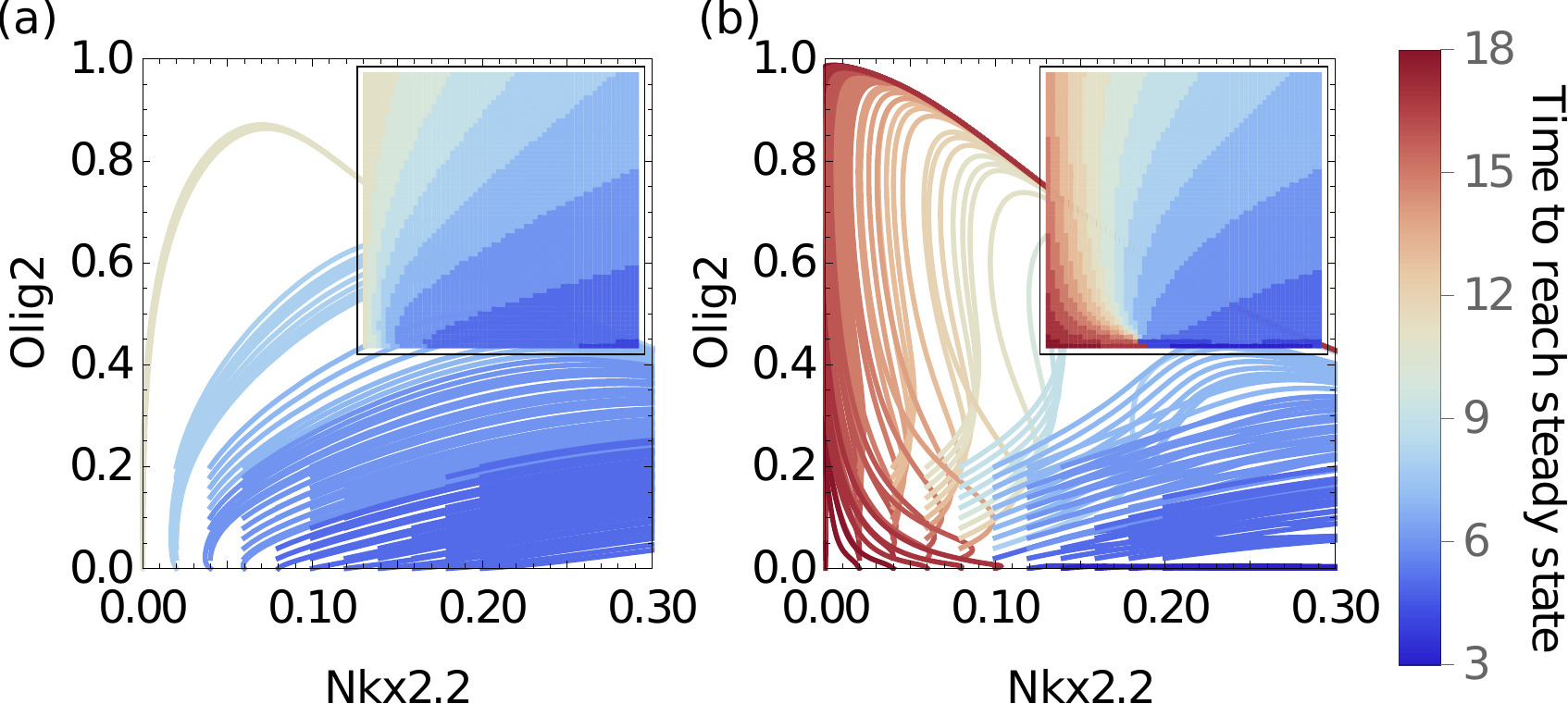}
\caption{Trajectories from multiple initial conditions starting near zero for Olig2 and Nkx2.2 in the p3 region ($p=0.1$). The comparison emphasises the robustness provided by the memory functions (ZMs, b) in comparison to a memoryless system (QSS, a). All trajectories reach the same attractor with high Nkx2.2 and low Olig2. The \zms system behaves almost identically to the full system, with transient increases in Olig2 and large delays in reaching the steady state from a variety of initial conditions. Colour indicates time taken to reach final steady state, quantified as Nkx2.2 deviating less than $1\%$ from its asymptotic concentration. Insets have same axes as main plots and show time to final steady state from a given initial position. \label{fig:strms}}
\end{figure}

Finally, we examined whether the effect of the two memory functions (one for Olig2 and one for Nkx2.2) that reflect the influence of the bulk is to increase the robustness of the system to initial conditions. For the system with memory, the delay in Nkx2.2 expression is present for multiple initial conditions, with trajectories crossing in a way that would be impossible to reproduce with a memoryless system (Fig.~\ref{fig:strms}). From Fig.~\ref{fig:mMem}a, the memory has two dominant channels reacting to changes in Nkx2.2 and Olig2, respectively.
This ensures that if even one of the subnetwork species levels drifts away from zero, the memory 
pushes the path back into the ``correct'' direction. In the case of the memoryless system, the already short transient observed in Fig.~\ref{fig:PONI_transient} disappears completely as soon as the initial conditions are no longer zero for Nkx2.2 (Fig.~\ref{fig:strms}a). In general a transient is difficult to achieve in a 2D memoryless system where a very specific function would have to repress Nkx2.2 at both low and medium-high levels of Olig2. We find that the memory generated by the combination of Pax6 and Irx3 provides robustness to changes in initial condition as the memory leads to low levels of Nkx2.2 during the initial phases of the transient (Fig.~\ref{fig:strms}b).

\section*{Discussion}
We have developed a version of the Zwanzig-Mori formalism, building on the work of Chorin and colleagues \cite{Chorin2000}, to obtain closed-form memory functions that can be used to reduce the dimensionality of a dynamical system far from equilibrium.
Our method applies to general first order nonlinear differential equations of arbitrary dimensionality. Systems with higher order time derivatives can in principle also be treated by introducing auxiliary variables for the lower order time derivatives in the conventional way, though we have not explored examples of this type.
The method is constrained only in that the bulk must have a unique steady state given the subnetwork state; we refer to this as the quasi-steady state (QSS). This restriction can be viewed as analogous to the assumption in the standard Zwanzig-Mori projection that, due to detailed balance, the bulk always reaches a Gibbs distribution for a given subnetwork state.

Similar to the method of averaging and other dimensionality reduction approaches \cite{Bronstein2018,Thomas2012,Sanders2007}, our method contains a timescale separation at its base, where the bulk is assumed to reach its QSS arbitrarily quickly, but in contrast it then systematically finds the memory terms that correct for this assumption.

%makes use of a timescale separation, our method then obtains systematic corrections around the limiting timescale-separated case.

In our examples the bulk parts of the network were usually chosen as small because, given the nonlinear cross-repressive  nature of the biological interactions that we study, multiple bulk steady states might otherwise result. However, there is no general restriction that the bulk has to be small per se or smaller than the subnetwork, as long as the assumption of a unique bulk steady state is met. For example, in the neural tube network, one could choose a single subnetwork species (Nkx2.2) and three bulk species. 

We have demonstrated the accuracy of the approach in capturing emergent dynamics such as non-trivial transients and sustained or damped oscillations. By construction the method can capture multistability and we have shown its accuracy in predicting the basins of attraction that, in a biological context, delineate cell fate decisions. We subsequently applied the method to a biologically relevant system, the neural tube patterning network \cite{Cohen2014}. The reduced model captures the nontrivial dynamics of the original system through its non-monotonic transients. In addition it provides a novel understanding of the cause of such transients in gene expression and suggests that memory effects, stored in independent bulk nodes, provide robustness to initial conditions.

As more information is acquired about complex systems, methods like ours that coarse grain the elements of a system but conserve the dynamics will be crucial to provide an understanding of such systems. We have demonstrated the generality of this method and its flexibility and applicability to dynamical systems. In the biological context, the approach holds promise e.g.\ for  more complex networks that incorporate signalling and gene regulation dynamics, where it could be applied to distinguish the impact of these two effects onto the macroscopic behaviour in e.g.\ tissue patterning.

In the present work we have not explored the role of the random force, which captures effects of the bulk initially being away from QSS. However, our method gives a closed-form approximation for this term in the projected equations, which is simply the function $F_s$ given by (\ref{memExp},\ref{fbs_sol}). This should make it straightforward to explore random force effects in the future, e.g.\ to assess the importance of changes in bulk initial conditions in the spirit of our previous work \cite{Herrera-Delgado2018}.  

A further interesting avenue of research would be the development of methods to capture the dynamics of systems even when parts of the network are unknown \cite{Xing2011,Mukhopadhyay2018}. This approach could then be coupled with identifying the contribution of specific memory channels to the observed dynamics.

% \section*{Materials and methods}
% The equations and corresponding parameters are either stated or given in the original articles cited.

% \section*{Author contributions}
% The manuscript was conceptualised, written and revised by EHD, JB and PS. EHD and PS performed all mathematical calculations.

% \section*{Competing interests}
% The authors declare no competing interests.

\section*{Acknowledgements}
We thank Matthias Fuchs and Alberto Pezzotta for helpful discussions and comments on the manuscript. This work was supported by: the Francis Crick Institute, which receives its core funding from Cancer Research UK (FC001051), the UK Medical Research Council (FC001051), and Wellcome (FC001051); funding from Wellcome [WT098325MA and WT098326MA]; the European Research Council under European Union (EU) Horizon 2020 research and innovation program grant 742138.

% Bibliography
\bibliography{library}

\begin{thebibliography}{10}

\bibitem{Davidson2010}
E.~H. Davidson, ``{Emerging properties of animal gene regulatory networks},''
  {\em Nature}, vol.~468, no.~7326, pp.~911--920, 2010.

\bibitem{Snider2015}
J.~Snider, M.~Kotlyar, P.~Saraon, Z.~Yao, I.~Jurisica, and I.~Stagljar,
  ``{Fundamentals of protein interaction network mapping},'' {\em Molecular
  Systems Biology}, vol.~11, no.~12, p.~848, 2015.

\bibitem{Rega2005}
G.~Rega and H.~Troger, ``{Dimension Reduction of Dynamical Systems: Methods,
  Models, Applications},'' {\em Nonlinear Dynamics}, vol.~41, pp.~1--15, aug
  2005.

\bibitem{Schnoerr2017}
D.~Schnoerr, G.~Sanguinetti, and R.~Grima, ``{Approximation and inference
  methods for stochastic biochemical kinetics - A tutorial review},'' {\em
  Journal of Physics A: Mathematical and Theoretical}, vol.~50, p.~093001, mar
  2017.

\bibitem{Bronstein2018}
L.~Bronstein and H.~Koeppl, ``{Marginal process framework: A model reduction
  tool for Markov jump processes},'' {\em Physical Review E}, vol.~97, no.~6,
  pp.~1--14, 2018.

\bibitem{Nakajima1958}
S.~Nakajima, ``{On Quantum Theory of Transport Phenomena},'' {\em Progress of
  Theoretical Physics}, vol.~20, pp.~948--959, dec 1958.

\bibitem{Zwanzig1961}
R.~Zwanzig, ``{Memory effects in irreversible thermodynamics},'' {\em Physical
  Review}, vol.~124, pp.~983--992, nov 1961.

\bibitem{Mori1965}
H.~Mori, ``{Transport, Collective Motion, and Brownian Motion},'' {\em Progress
  of Theoretical Physics}, vol.~33, pp.~423--455, mar 1965.

\bibitem{Kawasaki1973}
K.~Kawasaki and J.~D. Gunton, ``{Theory of nonlinear transport processes:
  Nonlinear shear viscosity and normal stress effects},'' {\em Physical Review
  A}, vol.~8, no.~4, pp.~2048--2064, 1973.

\bibitem{Zwanzig2001}
R.~Zwanzig, {\em {Nonequilibrium Statistical Mechanics}}.
\newblock Oxford University Press, 2001.

\bibitem{Chorin2000}
A.~J. Chorin, O.~H. Hald, and R.~Kupferman, ``{Optimal prediction and the
  Mori-Zwanzig representation of irreversible processes},'' {\em Proceedings of
  the National Academy of Sciences}, vol.~97, no.~7, pp.~2968--2973, 2000.

\bibitem{Chorin2002}
A.~J. Chorin, O.~H. Hald, and R.~Kupferman, ``{Optimal prediction with
  memory},'' {\em Physica D: Nonlinear Phenomena}, vol.~166, no.~3-4,
  pp.~239--257, 2002.

\bibitem{Weinan2008}
E.~Weinan, T.~Li, and E.~Vanden-Eijnden, ``{Optimal partition and effective
  dynamics of complex networks},'' {\em Proceedings of the National Academy of
  Sciences of the United States of America}, vol.~105, pp.~7907--7912, jun
  2008.

\bibitem{Chorin2006a}
A.~J. Chorin and P.~Stinis, ``{Problem reduction, renormalization, and
  memory},'' {\em Communications in Applied Mathematics and Computational
  Science}, vol.~1, pp.~1--27, dec 2006.

\bibitem{Stinis2006}
P.~Stinis, ``{A comparative study of two stochastic mode reduction methods},''
  {\em Physica D: Nonlinear Phenomena}, vol.~213, pp.~197--213, jan 2006.

\bibitem{Beck2009}
C.~L. Beck, S.~Lall, T.~Liang, and M.~West, ``{Model reduction, optimal
  prediction, and the Mori-Zwanzig representation of Markov chains},'' {\em
  Proceedings of the IEEE Conference on Decision and Control}, pp.~3282--3287,
  2009.

\bibitem{Thomas2012}
P.~Thomas, R.~Grima, and A.~V. Straube, ``{Rigorous elimination of fast
  stochastic variables from the linear noise approximation using projection
  operators},'' {\em Physical Review E - Statistical, Nonlinear, and Soft
  Matter Physics}, vol.~86, p.~041110, oct 2012.

\bibitem{Gouasmi2017}
A.~Gouasmi, E.~J. Parish, and K.~Duraisamy, ``{A priori estimation of memory
  effects in reduced-order models of nonlinear systems using the Mori-Zwanzig
  formalism},'' {\em Proceedings of the Royal Society A: Mathematical, Physical
  and Engineering Sciences}, vol.~473, no.~2205, p.~20170385, 2017.

\bibitem{Xing2011}
J.~Xing and K.~S. Kim, ``{Application of the projection operator formalism to
  non-Hamiltonian dynamics},'' {\em Journal of Chemical Physics}, vol.~134,
  no.~4, 2011.

\bibitem{Mukhopadhyay2018}
A.~Mukhopadhyay and J.~Xing, ``{A numerical procedure for model reduction using
  the generalized Langevin equation formalism},'' {\em arXiv}, apr 2013.

\bibitem{Rubin2014}
K.~J. Rubin, K.~Lawler, P.~Sollich, and T.~Ng, ``{Memory effects in biochemical
  networks as the natural counterpart of extrinsic noise},'' {\em Journal of
  Theoretical Biology}, vol.~357, pp.~245--267, sep 2014.

\bibitem{Bravi2017a}
B.~Bravi and P.~Sollich, ``{Inference for dynamics of continuous variables: The
  extended Plefka expansion with hidden nodes},'' {\em Journal of Statistical
  Mechanics: Theory and Experiment}, vol.~2017, p.~063404, mar 2017.

\bibitem{Kang2019}
H.~W. Kang, W.~R. KhudaBukhsh, H.~Koeppl, and G.~A. Rempa{\l}a,
  ``{Quasi-Steady-State Approximations Derived from the Stochastic Model of
  Enzyme Kinetics},'' {\em Bulletin of Mathematical Biology}, vol.~81, no.~5,
  pp.~1303--1336, 2019.

\bibitem{Cohen2014}
M.~Cohen, K.~M. Page, R.~Perez-Carrasco, C.~P. Barnes, and J.~Briscoe, ``{A
  theoretical framework for the regulation of Shh morphogen-controlled gene
  expression},'' {\em Development}, vol.~141, pp.~3868--3878, oct 2014.

\bibitem{Kawasaki1970}
K.~Kawasaki, ``{Kinetic Equations and Time Correlation of Critical
  Fluctuations},'' {\em Annals of Physics}, vol.~61, pp.~1--56, 1970.

\bibitem{Gotze1991}
W.~G{\"{o}}tze, ``{Aspects of structural glass transitions},'' {\em Liquids,
  Freezing and Glass Transition. Les Houches Summer Schools of Theoretical
  Physics Session LI}, pp.~287--503, 1991.

\bibitem{Herrera-Delgado2018}
E.~Herrera-Delgado, R.~Perez-Carrasco, J.~Briscoe, and P.~Sollich, ``{Memory
  functions reveal structural properties of gene regulatory networks},'' {\em
  PLoS Computational Biology}, vol.~14, no.~2, pp.~1--25, 2018.

\bibitem{Angeli2004}
D.~Angeli, J.~E. Ferrell, and E.~D. Sontag, ``{Detection of multistability,
  bifurcations, and hysteresis in a large class of biological positive-feedback
  systems},'' {\em Proceedings of the National Academy of Sciences}, vol.~101,
  no.~7, pp.~1822--1827, 2004.

\bibitem{Page2018}
K.~M. Page and R.~Perez-Carrasco, ``{Degradation rate uniformity determines
  success of oscillations in repressive feedback regulatory networks},'' {\em
  Journal of the Royal Society Interface}, vol.~15, no.~142, 2018.

\bibitem{Wilhelm2009}
T.~Wilhelm, ``{The smallest chemical reaction system with bistability},'' {\em
  BMC Systems Biology}, vol.~3, no.~1, p.~90, 2009.

\bibitem{Kondepudi1998}
D.~Kondepudi and I.~Prigogine, {\em {Modern Thermodynamics: From Heat Engines
  to Dissipative Structures}}.
\newblock John Wiley and Sons, 1998.

\bibitem{Caldeira1983}
A.~O. Caldeira and A.~J. Leggett, ``{Quantum tunnelling in a dissipative
  system},'' {\em Annals of Physics}, vol.~149, no.~2, pp.~374--456, 1983.

\bibitem{Zwanzig1973}
R.~Zwanzig, ``{Nonlinear generalized Langevin equations},'' {\em Journal of
  Statistical Physics}, vol.~9, no.~3, pp.~215--220, 1973.

\bibitem{Enver2009}
T.~Enver, M.~Pera, C.~Peterson, and P.~W. Andrews, ``{Stem Cell States, Fates,
  and the Rules of Attraction},'' {\em Cell Stem Cell}, vol.~4, pp.~387--397,
  may 2009.

\bibitem{Elowitz2000}
M.~B. Elowitz and S.~Leibier, ``{A synthetic oscillatory network of
  transcriptional regulators},'' {\em Nature}, vol.~403, no.~6767,
  pp.~335--338, 2000.

\bibitem{Strelkowa2010}
N.~Strelkowa and M.~Barahona, ``{Switchable genetic oscillator operating in
  quasi-stable mode},'' {\em Journal of The Royal Society Interface}, vol.~7,
  pp.~1071--1082, jul 2010.

\bibitem{Sagner2019}
A.~Sagner and J.~Briscoe, ``{Establishing neuronal diversity in the spinal
  cord: a time and a place.},'' {\em Development (Cambridge, England)},
  vol.~146, no.~22, 2019.

\bibitem{Exelby2019}
K.~Exelby, E.~Herrera-Delgado, L.~{Garcia Perez}, R.~Perez-Carrasco, A.~Sagner,
  V.~Metzis, P.~Sollich, and J.~Briscoe, ``{Precision of Tissue Patterning is
  Controlled by Dynamical Properties of Gene Regulatory Networks},'' {\em
  bioRxiv}, jul 2019.

\bibitem{Dessaud2007}
E.~Dessaud, L.~L. Yang, K.~Hill, B.~Cox, F.~Ulloa, A.~Ribeiro, A.~Mynett, B.~G.
  Novitch, and J.~Briscoe, ``{Interpretation of the sonic hedgehog morphogen
  gradient by a temporal adaptation mechanism.},'' {\em Nature}, vol.~450,
  pp.~717--20, nov 2007.

\bibitem{Sanders2007}
J.~A. Sanders, F.~Verhulst, and J.~Murdock, {\em {Averaging methods in
  nonlinear dynamical systems}}, vol.~59.
\newblock New York, NY: Springer, 2007.

\bibitem{John1978}
F.~John, {\em {Partial Differential Equations}}, vol.~1 of {\em Applied
  Mathematical Sciences}.
\newblock New York, NY: Springer US, 1978.

\end{thebibliography}

\clearpage
\onecolumn
\appendix
\setcounter{page}{1}
\renewcommand{\thefigure}{S\arabic{figure}}
\renewcommand{\theequation}{S\arabic{equation}}
\renewcommand{\thepage}{S\arabic{page}}
\setcounter{equation}{0}
\setcounter{figure}{0}
\setcounter{section}{0}

\section{Expansion around QSS \label{supp:expansion}}
In this appendix we detail how we derive the memory evolution over time. We make use of our expansion (\ref{memExp}) to find for the first term on the RHS of (\ref{memTime})
% we obtain the following expressions for $L\F$:
    \begin{align}
LF_s
%&+\sum_{s'}R_{s'}(\xs,\xb{})\sum_b\frac{\partial}{\partial x_{s'}}(x_b-x_b^*
%%(\xs)
%    )\fbs{b}{\tau}\nonumber\\
%%LF_s
    =&\sum_b R_b(\xs,\xb{})\fbs{b}{\tau} 
    +\sum_{s'}R_{s'}(\xs,\xb{})\sum_b(x_b-x_b^*
%(\xs)
)\frac{\partial}{\partial x_{s'}}\fbs{b}{\tau}\label{memTime_unsubtracted}\\
    &-\sum_{s'}R_{s'}(\xs,\xb{})\sum_b\frac{\partial x_b^*}{\partial x_{s'}}
%(\xs)
\fbs{b}{\tau}\nonumber
    \end{align}
where the two last terms arise by differentiating the product
%(\partial/\partial x_{s'})
$(x_b-x_b^*)\fbs{b}{\tau}$ w.r.t.\ $x_{s'}$ and we have not written the $\xs$-dependence of $\xb{*}$ for brevity. The second term on the RHS of (\ref{memTime}) is the expectation of this, obtained by replacing $\xb{}$ by $\xb{*}$:
%We then obtain the following similar expression for $E[L\F|\xs]$ where we use $\xb{*}$ instead of $\xb{}$
    \begin{align}
    E[L F_s|\xs]=&\sum_b R_b(\xs,\xb{*})\fbs{b}{\tau}
    -\sum_{s'}R_{s'}(\xs,\xb{*})\sum_b\frac{\partial x_b^*}{\partial x_{s'}}
%(\xs)
\fbs{b}{\tau}\label{memTimeExp}
    \end{align}
Putting the two together gives for the time evolution (\ref{memTime}) of $F_s$
%The first term actually vanishes by the definition (\ref{fastB}) of $\xb{*}$ but makes the structure 
%It's worth mentioning that the first term in (\ref{memTimeExp}) equals 0 as it describes the time differential of the fast variables as defined in (\ref{fastB}), yet we keep it as it is convenient for linearisation. We can then apply these substitutions of $L\F$ and $E[L\F|\xs]$ into the RHS of (\ref{memTime}) to obtain:
    \begin{align}
    \frac{\partial F_s}{\partial \tau}=&\sum_b\left[R_b(\xs,\xb{})-R_b(\xs,\xb{*})\right]\fbs{b}{\tau}
    +\sum_{s'}R_{s'}(\xs,\xb{})\sum_b(x_b-x_b^*
%(\xs)
)\frac{\partial}{\partial x_{s'}}\fbs{b}{\tau}\label{FOriginal}\\
    &-\sum_{s'b}\left[R_{s'}(\xs,\xb{})-R_{s'}(\xs,\xb{*})\right]\frac{\partial x_b^*}{\partial x_{s'}}
%(\xs)
\fbs{b}{\tau}\nonumber
    \end{align}
    For consistency with (\ref{memExp}) we now linearise the square brackets again in $\xb{}-\xb{*}$.
In the second term we similarly replace $R_{s'}(\xs,\xb{})$ by $R_{s'}(\xs,\xb{*})$ as the remaining factor in this term is already linear.
Comparing then with the time derivative of the original linearised formula (\ref{memExp}) gives, after appropriate relabelling of indices, the equation the evolution of $f_{bs}$ in time (\ref{fdto}).

From the above expansion one sees that the matrix $l_{bb'}$ in (\ref{fdto}) takes the form 
\begin{align}
    \lbb=\frac{\partial R_{b'}}{\partial x_b}-\sum_{s'}\frac{\partial R_{s'}}{\partial x_b}\frac{\partial x_{b'}^*}{\partial x_{s'}}
    \end{align}
The form (\ref{l_def}) in the main text is the obtained by using the identity
	\begin{align}
	\frac{\partial x_{b'}^*}{\partial x_{s'}}&=-\sum_{b''}(\bm{J}^{-1})_{b'b''}\frac{\partial R_{b''}}{\partial x_{s'}}\label{Jac2}
    \end{align}
The latter can be obtained by differentiating (\ref{fastB}) w.r.t.\ $x_s$.

To obtain the actual memory function from (\ref{eq:mem}) is now straightforward as we have already worked out the required expectation in (\ref{memTimeExp}). The first term on the r.h.s.\ of (\ref{memTimeExp}), which we had previously kept to make the ensuing linearisation easier to see, actually vanishes because of (\ref{fastB}), yielding
%  We will then simply use the fact the the first term of (\ref{memTimeExp}) is equal to $0$ in order to obtain:
\begin{align}
M_s(\xs,\tau)
%E[LF_s(\cdot,\tau)|\xs]
& = -\sum_{s'b'}R_{s'}\frac{\partial x_{b'}^{*}}{\partial x_{s'}}\fbs{b'}{\tau}
\end{align}
Using again the identity (\ref{Jac2}) we obtain our main result (\ref{M_final}).

\section{Solution for f \label{supp:f}}

In this appendix we find the solution $f_{bs}(\xs,\tau)$ to the differential equation (\ref{fdto}). We start by restating the latter as
\begin{align}   	
\frac{\partial}{\partial \tau}f_{bs}
-\sum_{s'}v_{s'}(\xs)\frac{\partial}{\partial x_{s'}}f_{bs}
&=\sum_{b'}\lbb(\xs) f_{b's}\label{fdto2}
\end{align}
where we have used that the factor $R_{s'}$ in (\ref{fdto}) is just the effective drift $v_{s'}$ defined in (\ref{subTime}).
As the equation is linear in $f_{bs}$ and its derivatives it can be solved using the method of characteristics (see e.g.~\cite{John1978}). Calling the curve parameter for a characteristic $u$, the characteristic equations can be read off from (\ref{fdto2}) as
\begin{align}
\frac{d\tau}{du}&= {}1\label{char_u}\\
\frac{dx_s}{du}&= {}-v_s(\xs)\label{char_xs}\\
\frac{df_{bs}}{du}&= {}\sum_{b'}\lbb(\xs) f_{b's} \label{char_fbs}
\end{align}
Setting an arbitrary integration constant to zero, the first of these gives $\tau=u$. To solve (\ref{char_xs}) we call $\bm{\phi}_v$ the flow generated by $\bm{v}(\xs)$, which is defined as the solution of the differential equation
\begin{align}
\frac{\partial}{\partial \tau}\bm{\phi}_v(\xs,\tau) = \bm{v}(\bm{\phi}_v(\xs,\tau)), \qquad \bm{\phi}_v(\xs,0)=\xs
\end{align}
The solution of (\ref{char_xs}) is then
\begin{equation}
\xs(u)=\bm{\phi}_v(\xs_0,-u)
\label{char_xs_sol}
\end{equation}
where $\xs_0$ is the value at the beginning of the characteristic curve ($u=0$); the minus sign in the second argument of $\bm{\phi}_v$ reflects the ``backward in time'' propagation in (\ref{char_xs}). We note for later that, as a consequence of (\ref{char_xs_sol}), the solution values at $u_1$ and $u_2$ are related by
\begin{equation}
\xs(u_2) = \bm{\phi}_v(\xs(u_1),-u_2+u_1)
\label{u1u2}
\end{equation}
Finally, the solution of (\ref{char_fbs}) is
\begin{equation}
f_{bs}(u)=\sum_{b'}\bigl(e^{\int_0^{u}du' \,\bm{l}(\xs(u'))}\bigr)_{bb'}f^0_{b's}(\xs_0)
\end{equation}
using the initial condition (\ref{f0_def}) at $\tau=u=0$. From (\ref{char_fbs}) we see that the matrix exponential appearing here must be time-ordered, with earlier ``times'' $u'$ appearing to the right of later ones.

It now remains to express $f_{bs}(u)$ in terms of $\xs(u)$ and $\tau(u)=u$. We fix a $\hat{u}=\hat{\tau}$ and call $\hxs=\xs(\hat{u})$. Using (\ref{u1u2}) with $u_2=u'$ and $u_1=u$ then shows that the $\xs$-solution (\ref{char_xs_sol}) can be expressed in terms of $\hxs$ as
\begin{equation}
\xs(u')=\bm{\phi}_v(\hxs,\hat{\tau}-u')
\end{equation}
and in particular $\xs_0=\bm{\phi}_v(\hxs,\hat{\tau})$, so that
\begin{equation}
f_{bs}(\hxs,\hat{\tau})=\sum_{b'}\big(e^{\int_0^{\hat{\tau}}du' \,\bm{l}(\bm{\phi}_v(\hxs,\hat{\tau}-u'))}\big)_{bb'}f^0_{b's}(\bm{\phi}_v(\hxs,\hat{\tau}))
\end{equation}
Changing integration variable to $\tau'=\hat{\tau}-u'$ and dropping the hats then gives the solution (\ref{fbs_sol}) announced in the main text. Note that as $\tau'=\hat{\tau}-u'$, the time ordering of the matrix exponential 
\begin{equation}
\bm{E}(\tau) = \exp\left({\int_0^\tau d\tau'\,\bm{l}(\bm{\phi}_v(\bm{x}^\textrm{s},\tau'))}\right)
\end{equation}
is such that the earlier $\tau'$ are now on the {\em left}. The appropriate time ordered matrix exponential is defined formally via its Taylor series
\begin{equation}
\bm{E}(\tau) = \bm{1} + \sum_{n=1}^\infty
\int \prod_{i=1}^n d\tau_i\,\bm{l}(\bm{\phi}_v(\xs,\tau_1))\times 
\cdots\times \bm{l}(\bm{\phi}_v(\xs,\tau_n))
%\int_{\tau_1}^{\tau_2} d\tau'\, \bm{l}(\bm{\phi}_v(\xs,\tau')) \\ %&+ \int_{\tau_1}^{\tau_2} d\tau' \int_{\tau'}^{\tau_2} d\tau''\,
%\bm{l}(\bm{\phi}_v(\xs,\tau'))
%\bm{l}(\bm{\phi}_v(\xs,\tau''))  + \ldots \nonumber
%%\expl{\tau_1}{\tau_2}
\end{equation}
with $\bm{1}$ the identity matrix and the integration in the other terms running over the range $0<\tau_1<\ldots<\tau_n<\tau$.
	
\section{Mapping of self-consistent memory to differential equations}\label{supp:selfcon}

We show in this appendix how to map the subnetwork equations with self-consistent memory,
	\begin{align}  
	\frac{\partial}{\partial t}x_s&=v_s(\xs(t))+
\tilde{\mathcal{M}}_s(t)
    \end{align}  
to a set of differential equations. The self-consistent memory term $\tilde{\mathcal{M}}_s(t)$ is given by (\ref{M_selfconsistent})
\begin{align}
\tilde{\mathcal{M}}_s(t) =& \sum_{b''}\int_0^{t}dt'\sum_{b'} c_{b'}(\xs(t'))\left(e^{\int_{t'}^t dt''\bm{l}(\xs(t''))}\right)_{b'b''} f_{b''s}^0(\xs(t))
\label{M_selfconsistent_app}
\end{align}
so can be written as
\begin{align}
\tilde{\mathcal{M}}_s(t) =& \sum_{b} m_{b}(t)f_{bs}^0(\xs(t))
\end{align}
with
\begin{align}
m_b(t) = 
\int_0^{t}dt'\ \sum_{b'} c_{b'}(\xs(t'))\left(e^{\int_{t'}^t dt''\bm{l}(\xs(t''))}\right)_{b'b}
\label{mb_integral}
\end{align}
It is then straightforward to check that
\begin{align}
\frac{d}{dt} m_b(t) = c_b(\xs(t)) + \sum_{b'} m_{b'}(t)l_{b'b}(\xs(t))
\label{dm_dt}
\end{align}
where the second term arises from the $t$-dependence of the matrix exponential. The $m_b(t)$ can therefore be obtained numerically by integrating the differential equations (\ref{dm_dt}) together with
%The full system of differential equations is then
% We then integrate the subnetwork species over time and for this we use a dummy variable $\bm{m}=\bm{c}\bm{E}$ as follows:
the subnetwork equations with (self-consistent) memory
\begin{align}
	\frac{d}{dt}x_s(t)&=v_s{(\xs(t))}+\sum_{b'}m_{b'}(t)f^{0}_{sb'}(\xs(t))
%\\
%\frac{d}{dt}m_{b}(t)&=c_{b}(\xs(t))+\sum_{b'}m_{b'}(t)l_{b'b}(\xs(t))
%\nonumber\\
%	\xs(0)&=\text{arbitrary}\nonumber\\
%	\bm{m}(0)&=0\nonumber
\end{align}
The appropriate initial conditions for the auxiliary variables follow from (\ref{mb_integral}) as $m_b(0)=0$.

\section{Channel decomposition} \label{sSec:decomp}

We begin by writing the expression for the memory function explicitly, combining Eqs.~(\ref{f0_def}, \ref{fbs_sol}, \ref{M_final}, \ref{c_def}):
\begin{align}
M_s(\xs,\tau) &=
\sum_{b's'b''}(\bm{J}^{-1})_{b'b''}\frac{\partial R_{b''}}{\partial x_{s'}} R_{s'}
f_{b's}(\xs,\tau)\\
&=
\sum_{b's'b''}(\bm{J}^{-1})_{b'b''}\frac{\partial R_{b''}}{\partial x_{s'}} R_{s'}
%f_{b's}(\xs,\tau)\\
\sum_{c}E_{b'c}(\tau) 
	f_{cs}^0(\bm{\phi}_v(\xs,\tau))\\
&=
\sum_{b's'b''}(\bm{J}^{-1})_{b'b''}\frac{\partial R_{b''}}{\partial x_{s'}} R_{s'}
\sum_{c}E_{b'c}(\tau) 
\frac{\partial R_s}{\partial x_c}(\bm{\phi}_v(\xs,\tau))
\end{align}
where the first three factors are evaluated at $\xs$. We now swap index labels and group the sums into a more intuitive form:
\begin{align}
M_s(\xs,\tau)
&=
%\sum_{b's'b''b}
\sum_{s'}\sum_{bb'}\Biggl(\sum_{b''}
(\bm{J}^{-1})_{b''b'}(\xs) \frac{\partial R_{b'}}{\partial x_{s'}}(\xs)
R_{s'}(\xs)
E_{b''b}(\tau) 
\frac{\partial R_s}{\partial x_b}(\bm{\phi}_v(\xs,\tau))\Biggr)\label{sEq:reOrdD}
%\\
%&=
%\sum_{s'}\sum_{bb'}\left(\sum_{b''}
%\frac{\partial R_s}{\partial x_b}(\bm{\phi}_v(\xs,\tau)) E_{b''b}(\tau) 
%(\bm{J}^{-1})_{b''b'}(\xs)\frac{\partial R_{b'}}{\partial x_{s'}}(\xs)R_{s'}(\xs)\right)%\nonumber
\end{align}
As discussed in the main text, the expression up to before the exponential represents a change in the deviation of the bulk species concentration $x_{b''}$ from its QSS values over some small time interval, in response to changes in the subnetwork concentrations $x_{s'}$ (see also (\ref{R_b_tilde}) below).
%The interpretation of this is as follows. A small change of subnetwork species $s'$, $R_{s'}dt'$, in a small time interval $dt'$ in the past first feeds through to a change of bulk species $b'$. 
In the factor $\partial R_{b'}/\partial x_{s'}$ only those bulk species $b'$ contribute whose time evolution depends explicitly on the subnetwork species $s'$ driving the bulk time evolution via $R_{s'}$.  The $b'$ can then be interpreted as {\em outgoing channels} for the signal from $s'$. After propagation in the bulk network the signal returns via another bulk species. Here only bulk species $b$ contribute that appear explicitly in the time evolution of subnetwork species $s$ as indicated by the factors $\partial R_s/\partial x_b$. The $b$ can therefore be interpreted as {\em incoming channels}. Overall, we have memory effects from $s'$ onto $s$, via an outgoing channel ($s'$ to $b'$) and an incoming channel ($b$ to $s$). Consistent with this interpretation, the outgoing channel ``susceptibilities'' $\partial R_{b'}/\partial x_{s'}$ are evaluated for the past, i.e.\ ``sending'', state $\xs\equiv\xs(t')$ of the subnetwork. The incoming channel susceptibilities $\partial R_s/\partial x_b$, on the other hand, are evaluated at the current time $t$ as shown by the propagation via $\bm{\phi}_v$ across the time difference $\tau=t-t'$. Within the self-consistent approximation (\ref{dm_dt}), this propagation corresponds directly to evaluation at the current state $\xs(t)$.

\section{Self-consistent channel decomposition \label{sSec:selfConDec}}

The channel decomposition of Sec.~\ref{supp:selfcon} can also be applied to the self-consistent memory approximation, as we now outline. Writing out the self-consistent memory term (\ref{M_selfconsistent_app}) explicitly and reordering and relabelling terms as in (\ref{sEq:reOrdD}) gives 
%\begin{align}
%\tilde{\mathcal{M}}_s(t) =& \sum_{b''}\int_0^{t}dt'\ \sum_{b'} c_{b'}(\xs(t'))\left(e^{\int_{t'}^t dt''\bm{l}(\xs(t''))}\right)_{b'b''}\nonumber\\&\times f_{b''s}^0(\xs(t))
%\label{M_selfconsistent_app0}
%\end{align}
%Decomposing as above produces
\begin{align}
\tilde{\mathcal{M}}_s(t) =& \int_0^{t}dt'\ 
\sum_{s'}\sum_{bb'}\sum_{b''}
\frac{\partial R_s}{\partial x_b}(\xs(t)) 
\left(e^{\int_{t'}^t dt''\bm{l}(\xs(t''))}\right)_{b''b}
(\bm{J}^{-1})_{b''b'}(\xs(t'))
\frac{\partial R_{b'}}{\partial x_{s'}}(\xs(t')) R_{s'}(\xs(t'))\\
=& \sum_{s'}\sum_{bb'} 
\frac{\partial R_s}{\partial x_b}(\xs(t)) 
m_{sbb's'}(t)\nonumber
\end{align}
where
\begin{align}
m_{sbb's'}(t) = \int_0^{t}dt'\ 
\sum_{b''}
(\bm{J}^{-1})_{b''b'}(\xs(t'))
\frac{\partial R_{b'}}{\partial x_{s'}}(\xs(t'))
R_{s'}(\xs(t'))
\left(e^{\int_{t'}^t dt''\bm{l}(\xs(t''))}\right)_{b''b}
\end{align}
From this last representation it follows that the $m_{sbb's'}(t)$ vanish at $t=0$ and obey the differential equations
\begin{align}
\frac{d}{dt} m_{sbb's'}(t) = {}
(\bm{J}^{-1})_{bb'}(\xs(t)) 
\frac{\partial R_{b'}}{\partial x_{s'}}(\xs(t))
R_{s'}(\xs(t))
+ \sum_{b''}m_{sb''b's'}(t)l_{b''b}(\xs(t))
\label{channel_DEs}
\end{align}
The channel-decomposed memory can therefore also be calculated from differential equations. Of course one only needs to find the $m_{sbb's'}$ for combinations $(sb)$ and $(b's')$ where the corresponding channel susceptibilities are non-zero.

\section{Exactness of memory \label{supp:exact}}
We show that the self-consistent memory $m_b(t)$ is exact when both $R_s$ and $R_b$ contain at most linear terms in $\xb{}$. 
In such a case, the full system can be written as
\begin{align}
R_s =& \ v_s + \sum_{b'} \tilde{x}_{b'} f^0_{b's},\qquad
R_b =\sum_{b'} J_{bb'} \tilde{x}_{b'}
\end{align}
where $\tilde{x}_b = x_b - x_b^*(\xs)$ and the QSS value $\xb{*}(\xs)$ is an arbitrary function of $\xs$. We now want to show that the $\tilde{x}_b$ correspond \emph{exactly} to the $m_b$ from the self-consistent \zms method. To do this we work out their evolution in time:
%[R_s(\xs,\xb{*})-J_{bb'}x_b](\bm{J}^{-1})_{bb'}
\begin{align}
\frac{d}{dt}\tilde{x}_b= {}R_b-\sum_{s'}R_{s'}\frac{\partial x_b^*}{\partial x_{s'}}
= \sum_{b'}J_{bb'}\tilde{x}_{b'}
+\sum_{s'}\left(v_{s'} + \sum_{b'} \tilde{x}_{b'} f^0_{b's'}\right)\sum_{b''}(\bm{J}^{-1})_{bb''}\frac{\partial R_{b''}}{\partial x_{s'}}
\end{align}
By using that for an $\xb{}$-linear system as assumed here one has $f_{bs}^0=\partial R_s/\partial x_b$, the above can be rewritten as:
\begin{align}
\frac{d}{dt}\tilde{x}_b=&\sum_{b''}(\bm{J}^{-1})_{bb''}\frac{\partial R_{b''}}{\partial x_{s'}}v_{s'}
+\sum_{b'} \tilde{x}_{b'}\left(J_{bb'}+\sum_{s'b''}(\bm{J}^{-1})_{bb''}\frac{\partial R_{b''}}{\partial x_{s'}}\frac{\partial R_{s'}}{\partial x_{b'}}\right)
\end{align}
Using then the definitions (\ref{l_def}~\&~\ref{c_def}) we obtain an expression equivalent to (\ref{dm_dt})
\begin{align}
\frac{d}{dt}\tilde{x}_b=c_b+\sum_{b'}\tilde{x}_{b'}l_{b'b}
\end{align}
thus showing that $\tilde{x}_b=m_b$ when we start from the same initial condition $\tilde{x}_b=0$, \emph{i.e.}\ the bulk at QSS.

\begin{figure}[h!]
\centering
\includegraphics[scale=0.5]{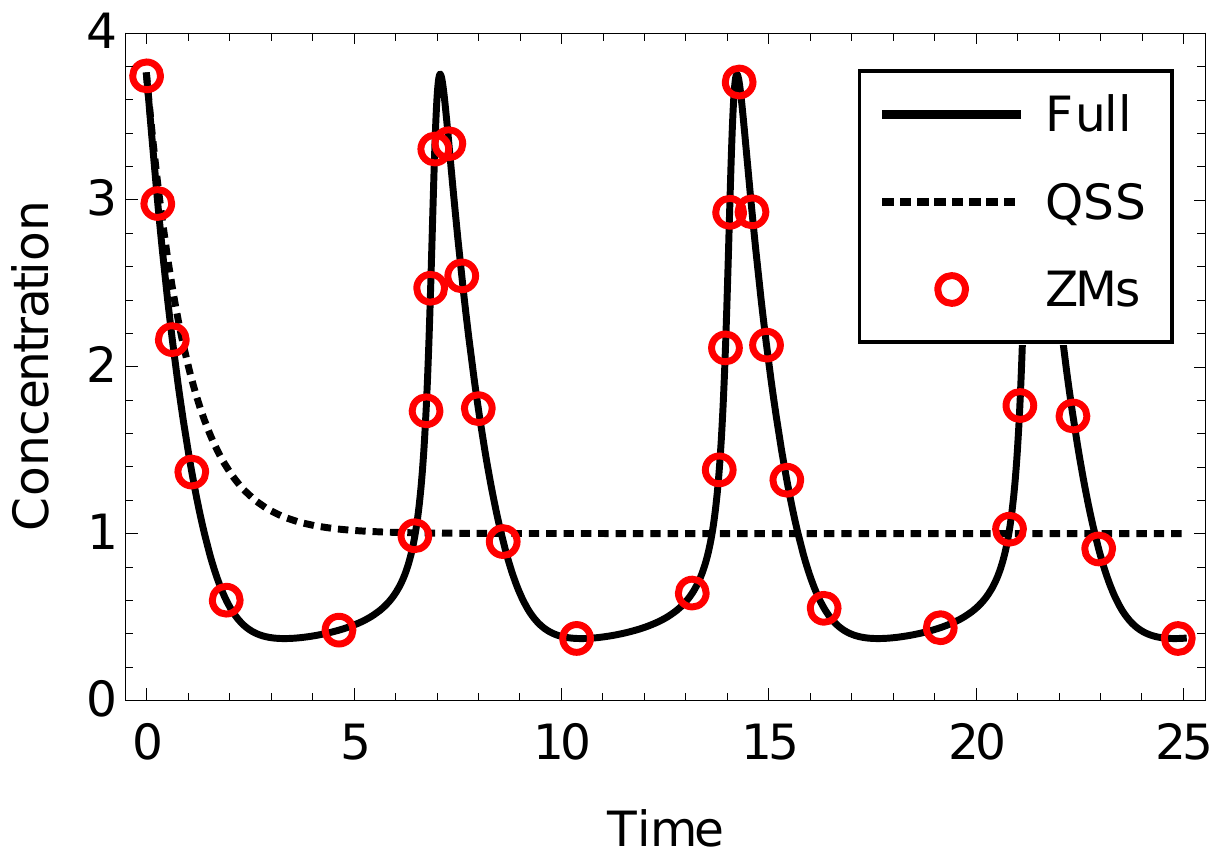}
\caption{
Oscillating Brusselator system as described in \cite{Kondepudi1998}; we retain the concentration of the first species $x_1$ in the subnetwork and place the second species in the bulk. Parameters for the sustained oscillatory regime in this example are $A=1$ and $B=3$ in the notation of \cite{Kondepudi1998}. The trajectory of the self-consistent projection (red dots) captures that of the original system (solid line) exactly as expected from the general proof in Supp.~\ref{supp:exact}, while the QSS (dotted line) fails qualitatively.\label{sFig:exactOsc}}
\end{figure}

\begin{figure}[h!]
\centering
\includegraphics[scale=0.5]{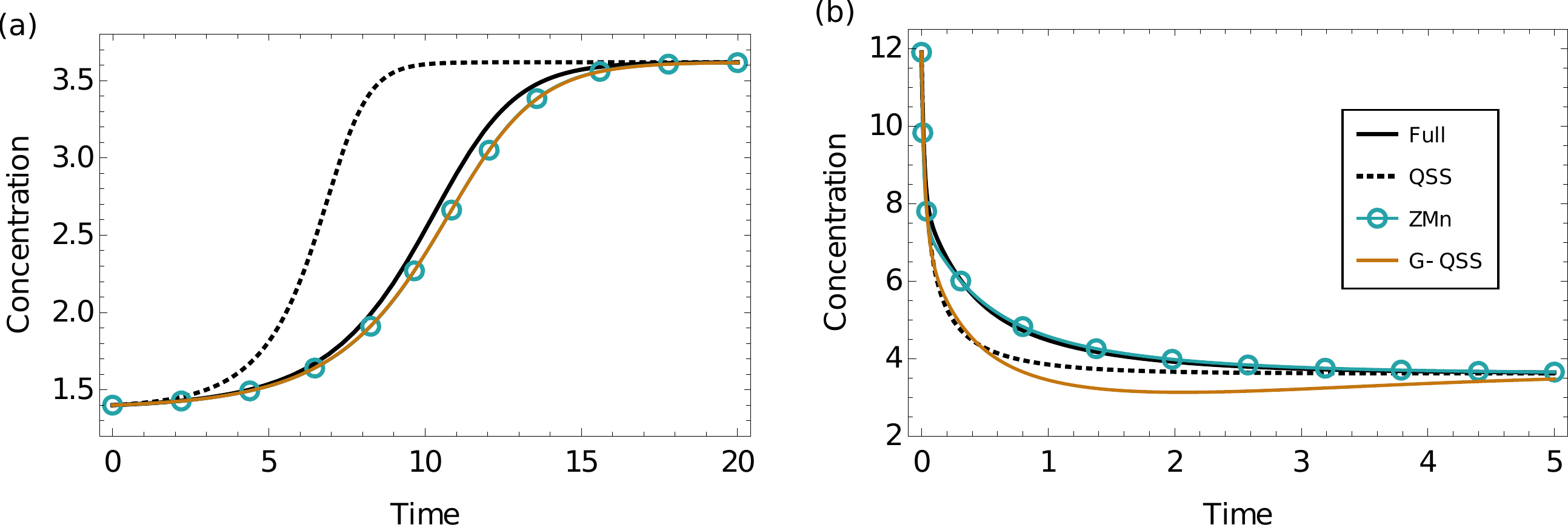}
\caption{
Minimal bistable system with linear dependence on the second species as described in \cite{Wilhelm2009}; we choose $x_1$ for the subnetwork and place the second species in the bulk. Parameters for the bistable regime in this example are: $k_1=10$, $k_2=1$, $k_3=2$ and $k_4=1$ \cite{Wilhelm2009}. 
The trajectory of the self-consistent \zms projection again captures that of the original system (solid line) exactly (see Supp.~\ref{supp:exact} \& Fig.~\ref{sFig:exactOsc}) so is not plotted. G-QSS represents the Gouasmi \emph{et al.} approach adapted to QSS projection (Supp.~\ref{sSec:GQSS}). (a) $x_1(0)=1.4$ At initial conditions near the fixed point both the \zmv method (cyan dots) and G-QSS (orange line) behave similarly and accurately capture the full dynamics. (b) $x_1(0)=11.9$ Further away from the final stable fixed point the G-QSS predictions become increasingly inaccurate while the \zmv method continues to provide a good approximation. \label{sFig:exactBis}}
\end{figure}

We test the above exactness statement on two different examples that have a linear dependence on a particular species but nonlinear dependences on other species: a minimal bistable system as described in \cite{Wilhelm2009}, and the ``Brusselator'', which is capable of achieving limit cycles \cite{Kondepudi1998}. As expected from the above derivation, the self-consistent memory captures the behaviour of both systems \textit{exactly} (Fig.~\ref{sFig:exactOsc} \& \ref{sFig:exactBis}). As further shown in Fig~\ref{sFig:exactBis}, the original nonlinear projection method \zmv is also accurate at capturing the dynamics though not necessarily exact. (We note that the memory functions of the Brusselator grow exponentially in a way that forces memory terms to cancel out to zero at the fixed point; this leads to numerical challenges that we do not pursue here.)

\section{Linear dynamics}
\label{supp:linear}

We discuss briefly the case of fully linear dynamics, where the dependence of $R_s$ and $R_b$ on {\em all} variables $\xs$ and $\xb{}$ (not just $\xb{}$ as in Supp.~\ref{supp:exact}) is only via constant and linear terms. Such a description can always be obtained by expanding linearly around a fixed point of the dynamics \cite{Rubin2014,Herrera-Delgado2018}. One then sees from (\ref{l_def}) and (\ref{f0_def}) that $l_{bb'}$ and $f^0_{bs}$ are both constant, \emph{i.e.}\ independent of $\xs$. Accordingly (compare (\ref{fbs_sol},~\ref{M_final}) and (\ref{M_selfconsistent})) also the memory functions of the \zmv and \zms projections become identical, and the corresponding channel decompositions are also the same.

\begin{figure}[h]
\centering
\includegraphics[scale=0.5]{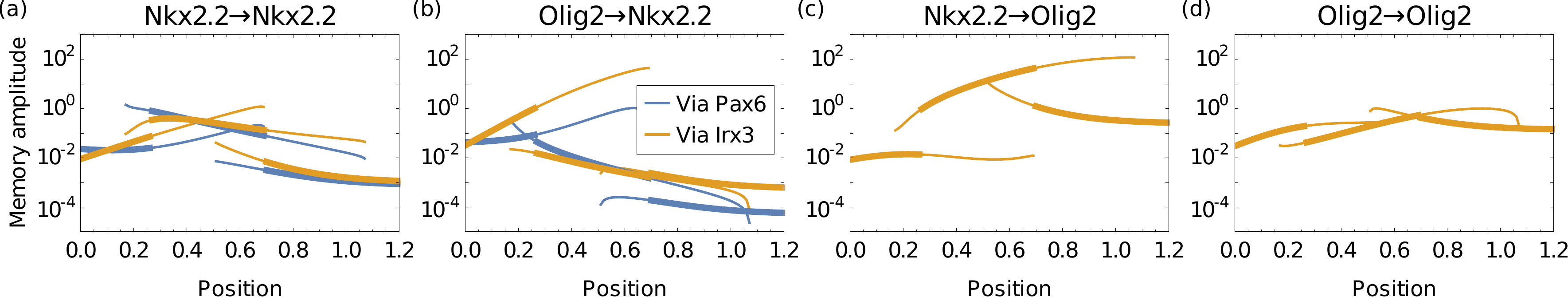}
\caption[Linear memory amplitude decomposition.]{
Amplitudes (value at time difference $\tau=0$) of nonlinear memory functions linearised around fixed points, for comparison with the approach of~\cite{Herrera-Delgado2018}, where memory functions expanded around fixed points were calculated directly. (a) Amplitude of memory of Nkx2.2 to itself along the neural tube. There are multiple lines as the analysis was performed at all possible stable steady states. The vertical axis is logarithmic to make the range of amplitudes easier to appreciate. Colours identify the memory amplitude contribution from the two possible bulk channels, via Irx3 and Pax6, respectively. Thick lines indicate physiological states, while thin lines indicate states that are not usually observed in vivo. (b) Amplitude of memory of Nkx2.2 to levels of Olig2, shown along the neural tube. The memory via Pax6 is for the most part below the memory via Irx3, in each pair of corresponding curves. (c, d) Amplitudes of memory of Olig2 to past Nkx2.2 (c) and to itself (d). No channel decomposition is performed as Olig2 receives memory only via the Irx3 channel. \label{fig:linChr}
}
\end{figure}

To illustrate the linearised dynamics approach, we perform a channel decomposition of the amplitude (value at $\tau=0$) of the linearised memory in the neural tube system as we did in \cite{Herrera-Delgado2018}, but now for the method derived in this study (Fig.~\ref{fig:linChr}). We find similar profiles to those found in \cite{Herrera-Delgado2018}. The results highlight the relative weakness of the memory from Olig2 into Nkx2.2 via Pax6, supporting the conclusions of \cite{Herrera-Delgado2018}. The method derived in this study is, however, significantly more powerful as it does not rely on an expansion near a steady state and gives access to the full memory and its channel decomposition as described in Supp.~\ref{sSec:decomp} \& \ref{sSec:selfConDec}.

\section{Comparison with alternative memory function approximation \label{sSec:Gouasmi}}

Gouasmi \textit{et al}~\cite{Gouasmi2017} propose an approximation for the memory function for the case where the projection (\ref{exp_def}) is defined not by setting the bulk coordinates to their $\xs$-dependent QSS values, but simply to zero:
\begin{equation}
E[g(\cdot) | \xs] = g(\xs,0)\label{exp_def_Gouasmi}
\end{equation}
The function $F_s(\xs,\xb{},\tau)$ still evolves according to (\ref{memTime}), which written out now reads
\begin{align}
\frac{\partial}{\partial\tau}F_s %(\xs,\xb{},\tau) 
=&LF_s(\xs,\xb{},\tau) - E[LF_s(\cdot,\tau)|\xs] \\
=&\sum_{s'}R_{s'}(\xs,\xb{})\frac{\partial F_s}{\partial x_{s'}}+\sum_b R_b(\xs,\xb{})\frac{\partial F_s}{\partial x_b}\\
&{}-\sum_{s'}R_{s'}(\xs,0)\frac{\partial F_s}{\partial x_{s'}}(\xs,0,\tau)
-\sum_bR_b(\xs,0)\frac{\partial F_s}{\partial x_b}(\xs,0,\tau)
\nonumber
\end{align}
where the very last factor is the $x_b$-derivative of $F$ evaluated at $\xb{}=0$. The approximation of~\cite{Gouasmi2017} amounts to ignoring the fact that the derivatives of $F_s$ are evaluated at a different point in the second line, which gives
\begin{align}
\frac{\partial}{\partial\tau}F_s=&\sum_{s'}\left[R_{s'}(\xs,\xb{})-R_{s'}(\xs,0)\right]\frac{\partial F_s}{\partial x_s}
+\sum_b \left[
R_b(\xs,\xb{})-R_b(\xs,0)\right]
\frac{\partial F_s}{\partial x_b}
\end{align}
This has the form of a Liouville equation as noticed in \cite{Gouasmi2017} and so its solution can be written as
%We write the memory function from as:
\newcommand{\psis}{\bm{\psi}^{\textrm{s}}}
\newcommand{\psib}{\bm{\psi}^{\textrm{b}}}
\begin{equation}
F_s(\xs,\xb{},\tau) = F_s(\psis(\xs,\xb{},\tau),
\psib(\xs,\xb{},\tau))
\end{equation}
where the components of the vector functions $\psis$ and $\psib$ evolve with $\tau$ according to
\begin{align}
\frac{\partial}{\partial\tau}\psi_s=&\ R_s(\psis,\psib)-R_s(\psis,0)
\label{psis_eq}\\
\frac{\partial}{\partial\tau}\psi_b=&\ R_b(\psis,\psib) -
R_b(\psis,0)
\label{psib_eq}
\end{align}
from the initial conditions
\begin{equation}
\psi_b(\xs,\xb{},0)=x_b, \qquad
\psi_s(\xs,\xb{},0)=x_s
\end{equation}
The function $F_s$ at $\tau=0$, which as before we write without a time argument, is given by the analogue of (\ref{F_time_zero}),
\begin{equation}
F_s(\xs,\xb{}) = R_s(\xs,\xb{})-R_s(\xs,0)
\label{F_time_zero_Gouasmi}
\end{equation}
The corresponding memory as defined in (\ref{eq:mem}) is
\begin{equation}
M^{\textrm{G}}_s(\xs,\tau)=\sum_{s'}R_{s'}\frac{\partial}{\partial x_{s'}}F_s(\bm{\psi}^{\textrm{s}},\bm{\psi}^{\textrm{b}})
+\sum_{b}R_{b}\frac{\partial}{\partial x_b}F(\bm{\psi}^{\textrm{s}},\bm{\psi}^{\textrm{b}})
\end{equation}
where all $R_{s'}$, $R_b$ and the derivatives are evaluated at $(\xs,0)$. Gouasmi \emph{et al} propose to find these derivatives numerically, but in fact a %n essentially
closed form expression can be obtained, as follows. Applying the chain rule gives
\begin{align}
M^{\textrm{G}}_s(\xs,\tau)=&\sum_{s's''}R_{s'}\frac{\partial F_s}{\partial \psi_{s''}}\frac{\partial\psi_{s''}}{\partial x_{s'}} +\sum_{s'b}R_{s'}\frac{\partial F_s}{\partial \psi_{b}}\frac{\partial\psi_{b}}{\partial x_{s'}}
+\sum_{s'b}R_{b}\frac{\partial F_s}{\partial\psi_{s'}}\frac{\partial \psi_{s'}}{\partial x_{b}}+\sum_{bb'}R_{b}\frac{\partial F_s}{\partial\psi_{b'}}\frac{\partial \psi_{b'}}{\partial x_{b}}\label{M_Gouasmi_chain}
\end{align}
Now note that in the final evaluation we always use $\xb{}=0$, which from the differential equations (\ref{psis_eq}, \ref{psib_eq}) implies $\psis=\xs$, $\psib=0$ for all $\tau$. Hence in particular $\psi_b$ is independent of $\xs$ and so $\partial \psi_b/\partial x_{s'}=0$. We also have $F_s(\xs,0)=0$ from (\ref{F_time_zero_Gouasmi}), which implies $\partial F_s/\partial \psi_{s'}=0$.
Only the last term from (\ref{M_Gouasmi_chain}) thus survives:
\begin{align}
M^{\textrm{G}}(\xs,\tau)=\sum_{bb'}R_{b}\frac{\partial F_s}{\partial\psi_{b'}}\frac{\partial \psi_{b'}}{\partial x_{b}}
\label{Gouasmi_vanilla_mem}
\end{align}
and it remains to find $\partial\psi_{b'}/\partial x_b$. By differentiating (\ref{psib_eq}) for $\partial \psi_{b'}/\partial \tau$ w.r.t.\ $x_b$ one finds
\begin{equation}
\frac{\partial}{\partial\tau}\frac{\partial \psi_{b'}}{\partial x_b} = \frac{\partial R_{b'}}{\partial \psi_{b''}}\frac{\partial \psi_{b''}}{\partial x_b}
\label{psi_response_eq}
\end{equation}
On the r.h.s.\ a similar term from the variation of $\psis$ vanishes as it would be proportional to
\begin{equation}
\frac{\partial R_{b'}}{\partial \psi_s}(\psis,\psib)
- \frac{\partial R_{b'}}{\partial \psi_s}(\psis,0)
\end{equation}
This difference is zero in the final evaluation at $\xb{}=0$ (which implies $\psib=0$). For the same reason the derivatives $\partial R_{b'}/\partial \psi_{b''} = \partial R_{b'}/\partial x_{b''}$ are evaluated at $(\xs,0)$ and constant in time $\tau$. Collecting these derivatives into a matrix $\bm{k}$  with elements $k_{b'b''}$ and using that $\partial \psi_b' / \partial x_b = \delta_{bb'}$ ($=1$ for $b=b'$ and $=0$ otherwise) at $\tau=0$ gives then as the explicit solution of (\ref{psi_response_eq})
\begin{equation}
\frac{\partial\psi_{b'}}{\partial x_b} = (e^{{\bm k}\tau})_{b'b}
\end{equation}
and inserting into (\ref{Gouasmi_vanilla_mem}) yields
\begin{align}
M^{\textrm{G}}(\xs,\tau)=\sum_{bb'}\frac{\partial R_s}{\partial x_{b'}}(e^{{\bm k}\tau})_{b'b}R_{b}
\label{memory_function_Gouasmi_vanilla}
\end{align}
where we have used that $\partial F_s/\partial \psi_{b'} = \partial F_s/\partial x_{b'} = \partial R_s/\partial x_{b'}$; this derivative and the factor $R_b$ are evaluated at $(\xs,0)$ in the approximation from \cite{Gouasmi2017} for
the memory function.

We do not show here how the above memory approximation performs in our test systems because the nature of the approach can lead to fixed points disappearing after projection or new fixed points appearing. We observed both of these effects in numerical evaluations for the bistable switch from \cite{Wilhelm2009}.

\section{Extending Gouasmi \emph{et al} approximation with QSS projection \label{sSec:GQSS}}

The Gouasmi \emph{et al.}\ approximation \cite{Gouasmi2017} for the memory function rests on projecting to $\xb{}=0$, but this is not generally an appropriate baseline for our case as it would correspond to setting all bulk concentrations to zero. However, we can adapt the approximation to the spirit of our work by changing coordinate system so that zero bulk coordinates correspond to the projection we consider throughout the paper, i.e.\ to QSS bulk concentrations. Explicitly, this variable transformation reads
\newcommand{\xt}{\tilde{x}}
\newcommand{\xvt}{\tilde{\bm{x}}}
\begin{equation}
\xt_s = x_s, \qquad \xt_b = x_b - x_b^*(\xs)
\end{equation}
because $\xt_b=0$ is then equivalent to $x_b = x_b^*(\xs)$. The time evolution of the new variables follows as
\begin{align}
\frac{d}{dt}\xt_s =&\ \tilde{R}_s(\xvt^{\rm s},\xvt^{\rm b}) \, =\, R_s(\xvt^{\rm s},\xb{*}+\xvt^{\rm b})\\
\frac{d}{dt}\xt_b =&\ \tilde{R}_b(\xvt^{\rm s},\xvt^{\rm b})\\
=&\ R_b(\xvt^{\rm s},\xb{*}+\xvt^{\rm b})
+
\sum_s\left[\sum_{b'}(\bm{J}^{-1})_{bb'}\frac{\partial R_{b'}}{\partial x_{s}}\right] R_s(\xvt^{\rm s},\xb{*}+\xvt^{\rm b})
\label{R_b_tilde}
\end{align}
where the factors enclosed in square brackets are the explicit expression for $-\partial x_b^*/\partial x_s$ and have to be evaluated at $\xb{*}$.

The Gouasmi memory function approximation, adapted for our QSS projection, is now given by (\ref{memory_function_Gouasmi_vanilla}) applied to the new variables $\xt_s$, $\xt_b$ and corresponding drift functions $\tilde{R}_s$, $\tilde{R}_b$. The last factor is $\tilde{R}_b(\xvt^{\rm s},0)$, which can be read off from (\ref{R_b_tilde}). The first term in (\ref{R_b_tilde}) vanishes as $R_b=0$ at QSS, while the remainder is seen to be precisely $c_b(\xs)$ from (\ref{c_def}). The matrix $\bm{k}$ in the new variables has elements
\begin{eqnarray}
k_{b'b} = \frac{\partial \tilde{R}_{b'}}{\partial \xt_b} = \frac{\partial R_{b'}}{\partial x_b}
%\nonumber\\
%&&{}
+\sum_s\left[\sum_{b'}(\bm{J}^{-1})_{bb'}\frac{\partial R_{b'}}{\partial x_{s}}\right]\frac{\partial R_s}{\partial x_b} = l_{bb'}
\end{eqnarray}
Note that the terms in square brackets are already just dependent on $\xs$, so do not contribute to the derivative. The remaining factor in the memory function is, again in the new variables,
\begin{equation}
\frac{\partial \tilde{R}_s}{\partial \xt_{b'}} 
 = \frac{\partial R_s}{\partial x_{b'}}
\end{equation}
so that overall
\begin{equation}
\tilde{M}^{\textrm{G}}_s(\xs,\tau) = \sum_{bb'} \frac{\partial R_s}{\partial x_{b'}}(e^{{\bm l}(\xs)\tau})_{bb'}c_{b}(\xs) = \sum_{b'} c_{b'}(\xs)\tilde f_{b's}(\xs,\tau)\label{M_g}
\end{equation}
with
\begin{equation}
\tilde f_{b's}(\xs,\tau) = \left(e^{{\bm l}(\xs)\tau}\right)_{b'b} f^0_{bs}(\xs)\label{f_G}
\end{equation}
where we have used the definition of $f^0_{bs}$ from (\ref{f0_def}). Comparing now (\ref{fbs_sol}) and (\ref{f_G}) shows that the memory approximation (\ref{M_g}), though derived here from rather different arguments, is quite similar to our expression (\ref{M_final}): the only difference is that the propagation from $\xs$ to $\bm{\phi}_v(\xs,\tau)$ is absent in $\tilde f_{b's}$, which is the analogue of our $f_{bs}$. We show that without the $\bm{\phi}_v$-propagation the method can still perform accurately in some situations but breaks down in other settings (Fig.~\ref{sFig:exactBis}).

\section{Neural tube model \label{sSec:NT}}
\newcommand{\xin}{x_{\rm in}}
We detail the model taken from \cite{Cohen2014} used to model ventral neural tube patterning.
The equations are as follows:
\begin{align}
\frac{d}{dt}x_\textrm{P}=&\frac{\alpha_{{\rm P}} w_{\textrm{P,p}}}{w_{\textrm{P,p}}+(1+k_{\textrm{PO}}x_{\textrm{O}})^2(1+k_{\textrm{PN}}x_{\textrm{N}})^2}-\beta_{\textrm{P}}x_\textrm{P}\\
\frac{d}{dt}x_\textrm{O}=&\frac{\alpha_{{\rm O}} w_{\textrm{O,p}}(1+k_{\textrm{O,in}}\xin)}{w_{\textrm{O,p}}(1+k_{\textrm{O,in}}\xin)+(1+k_{\textrm{OI}}x_{\textrm{I}})^2(1+k_{\textrm{ON}}x_{\textrm{N}})^2}-\beta_{\textrm{O}}x_\textrm{O}\\
\frac{d}{dt}x_\textrm{N}=&\frac{\alpha_{{\rm N}} w_{\textrm{N,p}}(1+k_{\textrm{N,in}}\xin)}{w_{\textrm{N,p}}(1+k_{\textrm{N,in}}\xin)+(1+k_{\textrm{NP}}x_{\textrm{P}})^2(1+k_{\textrm{NO}}x_{\textrm{O}})^2(1+k_{\textrm{NI}}x_{\textrm{I}})^2}-\beta_{\textrm{N}}x_\textrm{N}\\
\frac{d}{dt}x_\textrm{I}=&\frac{\alpha_{{\rm I}} w_{\textrm{I,p}}}{w_{\textrm{I,p}}+(1+k_{\textrm{IO}}x_{\textrm{O}})^2(1+k_{\textrm{IN}}x_{\textrm{N}})^2}-\beta_{\textrm{I}}x_\textrm{I}\label{eq:PONI}
\end{align}
The parameters and their meaning in a biological sense are detailed as follows:
\begin{center}
\hspace*{-0.7cm}
\begin{tabular}{| c | c | c |}
\hline
Name & Meaning & Value \\\hline
$\alpha_{\textrm{P}}$ & Pax6 production rate & 2 \\\hline
$\alpha_{\textrm{O}}$ & Olig2 production rate & 2\\\hline
$\alpha_{\textrm{N}}$ & Nkx2.2 production rate & 2 \\\hline
$\alpha_{\textrm{I}}$ & Irx3 production rate & 2 \\\hline
$\beta_{\textrm{P}}$ & Pax6 degradation rate & 2 \\\hline
$\beta_{\textrm{O}}$ & Olig2 degradation rate & 2 \\\hline
$\beta_{\textrm{N}}$ & Nkx2.2 degradation rate & 2 \\\hline
$\beta_{\textrm{I}}$ & Irx3 degradation rate & 2 \\\hline
$k_{\textrm{PO}}$ & Olig2 binding to Pax6 DNA & 1.9 \\\hline
$k_{\textrm{PN}}$ & Nkx2.2 binding to Pax6 DNA & 26.7 \\\hline
$k_{\textrm{ON}}$ & Nkx2.2 binding to Olig2 DNA & 60.6 \\\hline
$k_{\textrm{OI}}$ & Irx3 binding to Olig2 DNA & 28.4 \\\hline
$k_{\textrm{NP}}$ & Pax6 binding to Nkx2.2 DNA & 4.8 \\\hline
$k_{\textrm{NO}}$ & Olig2 binding to Nkx2.2 DNA & 27.1 \\\hline
$k_{\textrm{NI}}$ & Irx3 binding to Nkx2.2 DNA & 47.1 \\\hline
$k_{\textrm{IO}}$ & Olig2 binding to Irx3 DNA & 58.8 \\\hline
$k_{\textrm{IN}}$ & Nkx2.2 binding to Irx3 DNA & 76.2 \\\hline
$w_{\textrm{P,p}}$ & Polymerase binding to Pax6 DNA & 3.84 \\\hline
$w_{\textrm{O,p}}$ & Polymerase binding to Olig2 DNA & 2.01263 \\\hline
$w_{\textrm{N,p}}$ & Polymerase binding to Nkx2.2 DNA & 0.572324 \\\hline
$w_{\textrm{I,p}}$ & Polymerase binding to Irx3 DNA & 18.72 \\\hline
$k_{\textrm{O,in}}$ & Gli (Shh signal) binding to Olig2 DNA & 180 \\\hline
$k_{\textrm{N,in}}$ & Gli (Shh signal) binding to Nkx2.2 DNA & 373 \\\hline
\end{tabular}
\hspace*{-0.7cm}
\end{center} 
The signal input concentration $\xin$ is the gradient $e^{-s/0.15}$, which depends on the dorsal-ventral neural tube position $s$ ranging from 0 to 1 as in \cite{Cohen2014}.
\clearpage

\begin{figure}[h!]
\centering
\includegraphics[scale=0.5]{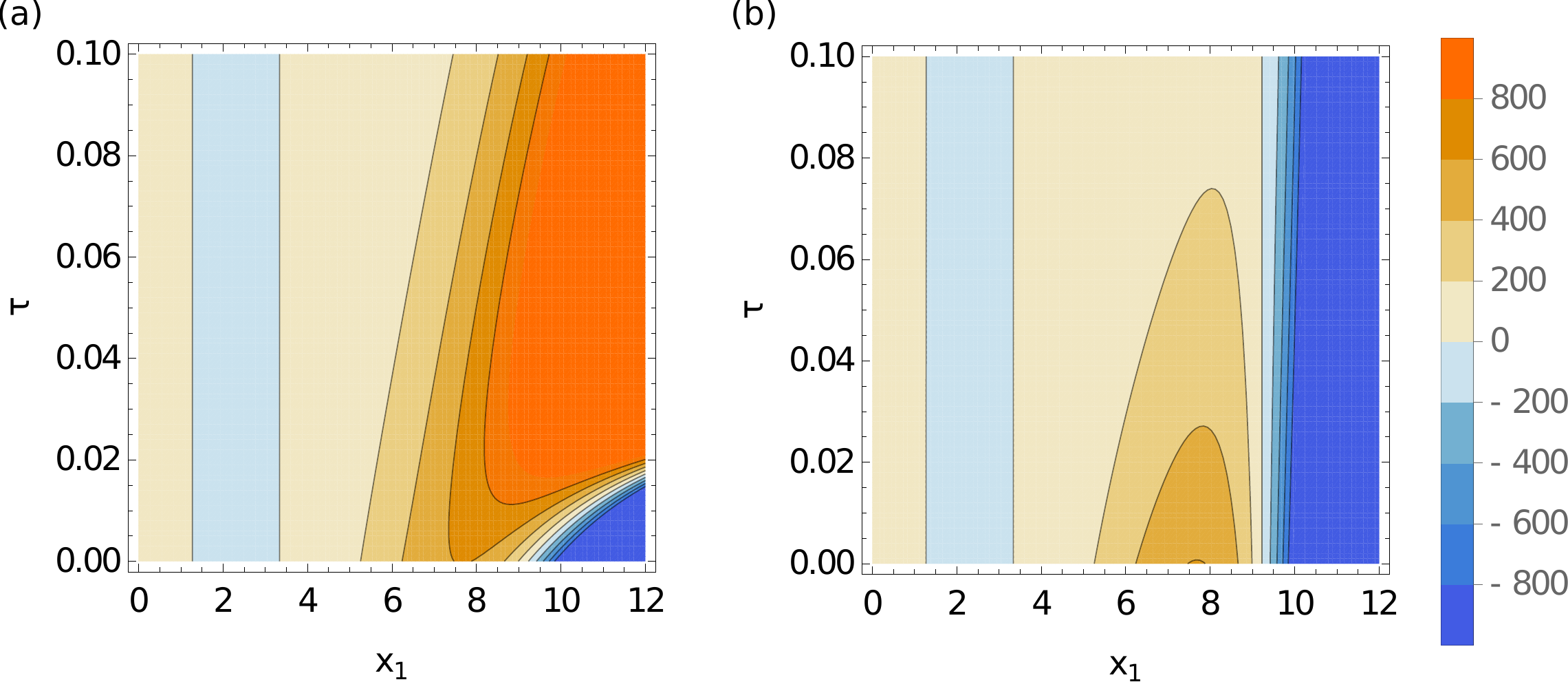}
\caption{Memory functions for the system detailed in \cite{Wilhelm2009} with the parameters chosen for Fig.~\ref{sFig:exactBis}. (a) Using the \zmv. (b) Using the method from Gouasmi \textit{et al.} as extended to QSS projection in Supp.~\ref{sSec:GQSS}. The $x$-axis shows the concentration of the subnetwork species $x_1$ while the $y$-axis indicates time difference $\tau$. By construction, the two memory function approximations predict the same value (scale bar to the right) at $\tau=0$, as they only differ in how they propagate the memory over time. At $\tau>0$ the memory functions are relatively similar for $x_1\in [0,6]$ but become progressively different as $x_1$ grows beyond this range; for $x_1 \geq 10$ the G-QSS method predicts a negative memory function for all $\tau$ that leads to its poor performance as
%$x_1$ grows larger the \zmv memory presents a sign change with $\tau$ that leads to the greater accuracy of the ZMn method
observed in Fig.~\ref{sFig:exactBis}. \label{sFig:bisMem}}
\end{figure}

\begin{figure}[h!]
\centering
\includegraphics[scale=0.5]{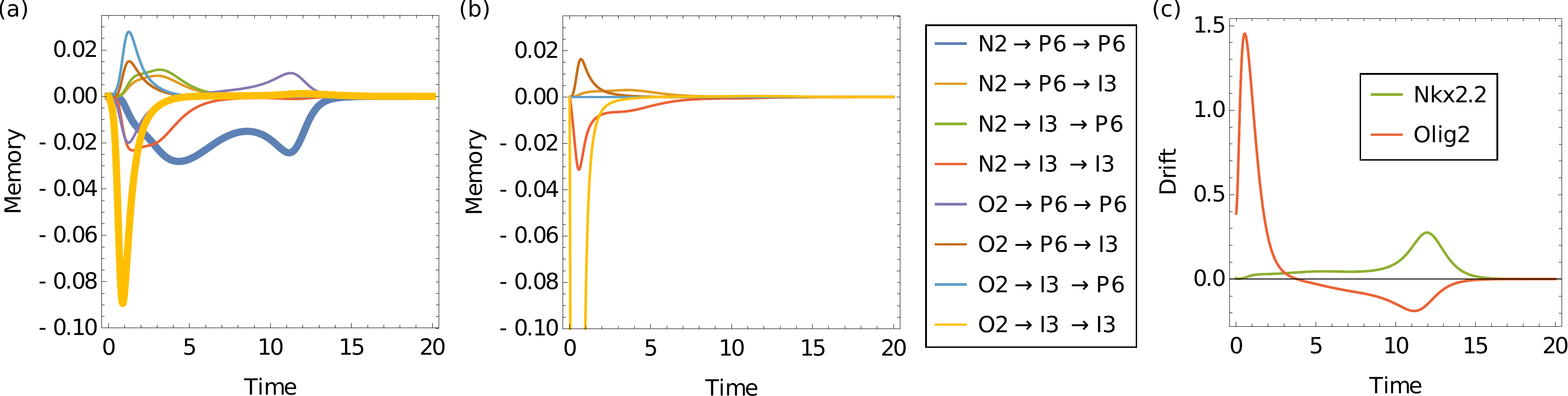}
\caption{Channel decomposition of the memory terms for the \zms trajectory shown in Fig.~\ref{fig:PONI_transient}. (a) and (b) show memory terms affecting Nkx2.2 and Olig2 respectively. Each colour represents a memory channel as indicated. The memory originates from a particular subnetwork species ``sending'' memory through a specific bulk species; the effect then propagates within the bulk and returns via a specific bulk species (see legend). The two most salient memory functions are: Nkx2.2 to Pax6 and then returning through Pax6 into Nkx2.2 (thick yellow line), and Olig2 to Irx3 and then returning through Pax6 into Nkx2.2 (thick blue line). (b) shows a large memory contribution acting on Olig2 via the channel through Irx3 (yellow line). However, in this case the drift for Olig2 (c) is so large that the relative effect of this memory channel remains nonetheless small.
 \label{sFig:decompP3}}
\end{figure}

\end{document}